\documentclass{singlecol-new}

\usepackage{natbib,stfloats}
\usepackage{mathrsfs}
\usepackage{float}
\usepackage{graphicx}
\usepackage{amsmath} 
\usepackage{amssymb}

\theoremstyle{TH}{
\newtheorem{lemma}{Lemma}
\newtheorem{theorem}[lemma]{Theorem}

}

\theoremstyle{THrm}{

}

\theoremstyle{THhit}{

}

\makeatletter
\def\theequation{\arabic{equation}}

\makeatother

\begin{document}%

\newcommand{\be}{\begin{equation}}
\newcommand{\ee}{\end{equation}}
\newcommand{\bea}{\begin{eqnarray}}
\newcommand{\eea}{\end{eqnarray}}
\newcommand{\bee}{\begin{eqnarray*}}
\newcommand{\eee}{\end{eqnarray*}}
\newcommand{\D}{${\bf{D~}}$}
\newcommand{\C}{${\bf{C~}}$}
\newcommand{\0}{${\bf 0}$}
\newcommand{\nii}{\noindent}
\newcommand{\nn}{\nonumber}

\setcounter{page}{1}









\title{On state dependent batch service queue with single and multiple vacation under Markovian arrival process}

%
%
%
%
%

%
%
%
%
%
%
%

\begin{abstract} An infinite buffer batch service vacation queue has been studied where service rate of the batch is dependent on the size of the batch and vacation rate is dependent on the queue size at vacation initiation epoch. The arrivals follow the Markovian arrival process (MAP). For service rule, general bulk service (GBS) rule is considered. The service time and vacation time both are considered to be generally distributed. Several joint distributions of interest are obtained using the bivariate vector generating function method and the supplementary variable technique (SVT). Numerical results are presented to show the behavior of the system performance to validate the analytical results.
\end{abstract}

\KEYWORD{MAP; Infinite buffer; GBS rule; Bivariate VGF; SVT.}

%
%
%

\maketitle

 \section{Introduction} Queues in which the server processes customers in groups, known as bulk (batch) service queues. In the study of bulk service queues, several kinds of bulk service rules are well known \emph{viz.,} fixed batch size bulk service rule, general bulk service (GBS) or ($a,b$) rule, random batch size bulk service rule, etc. In addition to its many applications in telecommunication, blood testing, etc., the GBS rule proposed by \cite{neuts1967general} is most familiar. Under the GBS rule, the server determines the thresholds $a$ (minimum) and $b$ (maximum), respectively, to serve the customers. Once service is completed, if the customers waiting in queue are between $a$ and $b$, all the customers are serviced. If the customers waiting in queue are more than $b$, it takes $b$ customers for the service, and if the customers waiting in queue are $<a$, then it falls into a dormant state until the queue size $a$ is reached. Apart from this, in the vacation queue proposed by \cite{levy1975utilization} the server may use the dormant period for some secondary work (\emph{viz.,} maintenance, promoting the revised policy of the company, etc.). In the single vacation (SV) queue, after a service, if the required number of customers is not present in the queue, then the server is sent on vacation. After returning from vacation, if the required number of customers is present in the queue, the customer is operated; otherwise, it stays dormant in the system. However, in the multiple vacation (MV) queue, the server takes repeated vacations until it finds the necessary number of customers in the system at the end of the vacation. Due to their direct applicability in transportation, telecommunication, production management, etc., bulk service queues with vacations have greatly impacted the researchers. For current work in bulk service queues with vacation, readers are invoked to see \cite{tamrakar2020steady,tamrakar2021study,gupta2020finite,nandy2021joint,deepa2022analysis} and the references therein. Remarkable survey carried out on vacation model by \cite{doshi1986queueing} and \cite{ke2010recent}. The literature that provides a key understanding on bulk service (vacation) queues readers are invited to see the books by \cite{chaudhry1983first} and \cite{medhi2002stochastic} (\cite{takagi1991queueing} and \cite{tian2006vacation}).

Bulk service non vacation queue where the arrival process is non-renewal  analyzed by few researchers, see \cite{BANERJEE2015138,pradhan2017analysis,pradhan2020stationary,singh2013computational,bank2020analytical,samanta2020analysis}.
\cite{BANERJEE2015138}, \cite{pradhan2017analysis}, and \cite{pradhan2020stationary}, respectively, considered $MAP/G_r^{(a,b)}/1/N$, $MAP/G_r^{(a,b)}/1$ and $BMAP/G_r^{(a,b)}/1$ queue and obtained the joint distribution of queue and number with server at various epochs. \cite{singh2013computational}, \cite{bank2020analytical} and \cite{samanta2020analysis}, respectively, analyzed $MAP/R^{(a,b)}/1$, $BMAP/G^{(a,Y)}/1$, and $BMAP/R^{(a,b)}/1$ queue for queue length distribution at different epochs. The authors of \cite{BANERJEE2015138, pradhan2017analysis, pradhan2020stationary, singh2013computational, bank2020analytical, samanta2020analysis} conducted their analysis using the SVT, embedded Markov chain technique (EMCT) and roots method.

In most of the bulk service queueing system with vacation authors have been considered Poisson/renewal arrival process, see, e.g.,  \cite{sikdar2016analysis,gupta2020finite,ayyappan2020analysis,nandy2021joint,chang2005performance,sikdar2008batch,haridass2012analysis,jeyakumar2014modelling,ayyappan2018analysis}, however, in most of the queues (for ex. telecommunication, computer network, etc.) the Poisson/renewal arrival process do not fit due to highly irregular traffic. A good representation for analyzing such bursty and correlated traffic is a non-renewal arrival process, i.e., the Markovian arrival process (MAP) proposed by \cite {lucantoni1990single}.  Some other input process is also included in MAP, \textit{viz}., Markov modulated Poisson process (MMPP), the phase (PH)-type renewal process, the interrupted
Poisson process (IPP), Poisson process. \cite{gupta1900finite} and \cite{sikdar2016analysis}, respectively, discussed $MAP/G^{(a,b)}/1/N$ queue with SV and $BMAP/G^Y/1/N$ queue with SV (MV), respectively, and obtained queue length distribution at various epoch using EMCT and SVT. In discrete-time set-up, \cite{nandy2021joint} considered discrete-time batch size dependent batch service queue with SV and MV. They carried out their analysis for queue and server content and queue length probabilities using the SVT and bivariate probability generating function technique.

  In most of the vacation queueing models, authors considered the length of vacation of the server as random and unaffected by the queue size at the vacation initiation epoch. The queueing model with vacation where the length of vacation depends on the queue size (length) at the vacation initiation epoch is known as the queue size dependent (QSD) vacation model. Such a QSD vacation models are analyzed by a few researchers, see \cite{shin1998bmap,banik2013analysis,thangaraj2017analysis,gupta2020finite,tamrakar2021study}. \cite{gupta2020finite} considered  $M/G_r^{(a,b)}/1/N$ queue with QSD SV (MV) and obtained the joint distribution of queue and server content and the joint distribution of queue content and vacation type using the SVT. Later, \cite{tamrakar2021study} analyzed the same model as \cite{gupta2020finite} for infinite buffer queue using the SVT and bivariate generating function approach. To the extent of the authors' knowledge, an infinite capacity queueing system $MAP/G_r^{(a,b)}/1$ queue with queue size dependent SV and MV has not been investigated previously in the literature. Further, we have considered that the queue size dependent vacation policy remarkably reduces congestion.

 In a computer network with highly irregular traffic, the proposed model is applicable. A desktop computer system connects to a local area network (LAN) via Ethernet (IEEE802.11h) link. Digital signals are transmitted over Ethernet in the form of a group (packet), with the transmission rate varying with the packet under transmission. Power utility (power utility increases with transmission rate) depends on transmission rate. The medium access control (MAC) handshake protocol is helpful in decreasing the average power utility. It is achieved by measuring the queue size (signals) waiting for transmission. After a transmission, if the number of signals lower than the previously established lower threshold, the handshake mechanism (vacation period) activates, and it depends on the queue size at the vacation initiation epoch.

 A model description of the considered model can be found in section \ref{sec2}. In Section \ref{sec3}, we have analyzed the model mathematically. Section \ref{sec4} presents some marginal distributions. Section \ref{sec5} contains the performance measures. Numerical results are presented  in Section \ref{sec6}, and for the conclusion, we present Section \ref{sec8}.

\section{Model}\label{sec2}
The single server with infinite capacity, batch size-dependent bulk service queue, and queue size-dependent vacation (single and multiple) is introduced. Markovian arrival process (MAP) governs the customer's arrival to the system, which is governed by the underlying Markov chain (UMC). In the UMC, there is a change from state $i$ to state $j$ ($1\leq i,j\leq m$). We set $d_{i,j}$ to be the transition rate with an arrival and $c_{i,j}$ to be the transition rate without an arrival. The $m\times m$ matrix $C=(c_{i,j})$ has non-negative off diagonal and negative diagonal members, whereas the $m\times m$ matrix $D=(d_{i,j})$ only has non-negative elements. The infinitesimal generator of the UMC is presented by the matrix ($C+D$). Assume that $\xi=( \xi_1, \xi_2,..., \xi_m)$ is a stationary probability vector such that $\xi($C+D$)=\textbf{0}$, $\xi \bf e$=1, where $\bf e$ is $m\times 1$ column matrix with each element 1 and \textbf{0} is $1\times m$ zero matrix. The fundamental arrival rate is determined by $\lambda= \xi D\bf e$. We suppose that $I$ refers to an identity matrix with an appropriate dimension. According to GBS rule, the consumers are serviced in batches (groups). The service time $(T_r)$ of a batch of size $r$ $(h\leq r\leq H)$ is generally distributed with probability density function (pdf) $s_r(t)$, distribution function (DF) $S_r(t)$, the Laplace-Satieties transform (LST) $\tilde{S}_r{(\theta)}$ and the mean service time $\frac{1}{\mu_r}=s_r=-\tilde{S}_r^{(1)}{(0)}$ $(h\leq r \leq H)$, where $\tilde{S}_r^{(1)}{(0)}$ is the derivative of $\tilde{S}_r{(\theta)}$ evaluated at $\theta$=0. When a service is finished, and if the server determines that the queue length is $l$ $(\geq h)$, it begins service in accordance with the GBS rule, i.e., it serves a batch of size min$(l,H)$, where $H$ is the server's maximum capacity. If the queue length is $k$ $(<h)$ after a service, the server begins the vacation, which has a random length and is dependent on the queue length $k$ $(0 \leq k \leq h-1)$.
 Let ${V}_k{(t)}$ \{$v_k{(t)}$\} [$\tilde{V}_k{(\theta)}$] be the DF  \{pdf\} [LST] of a typical vacation time $V_k$ $(0\leq k\leq h-1)$ which is generally distributed. The mean vacation time $\frac{1}{\nu_k}=x_k=-\tilde{V}_k^{(1)}{(0)}$ where $\tilde{V}_k^{(1)}{(0)}$ is the derivative of $\tilde{V}_k{(\theta)}$ at $\theta$=0. If the server finds at least $h$ waiting customers at the end of the vacation, it operates those customers in accordance with GBS rule; otherwise, it enters a state of dormancy until the queue length reaches the minimum threshold $h$, or it takes another vacation depending on the vacation policy being considered, namely either single vacation (SV) or multiple vacation (MV). The system's stability is ensured by the traffic intensity, which is $\frac{\lambda s_H}{H}(<1)$.
 Using the following definition of the variable $\epsilon$, we have examined single vacation (SV) and multiple vacation (MV) queues in this article in a unified manner.
 \begin{equation*}
 \epsilon=
 \begin{cases}
 1, &\mbox{for~ MV},\\
 0,&\mbox{for~ SV}.
 \end{cases}
 \end{equation*}

 \section{System analysis}\label{sec3} The following random variables at time $t$ are necessary for the mathematical analysis of our considered model.
  \begin{itemize}
\item $q(t)$ $\equiv$ Size of customers waiting in line (queue).
\item $S(t)$ $\equiv$ Customer size with the server when the server is busy.
\item $K_{v}(t)$ $\equiv$ Vacation type, when the server is on vacation.
\item $J(t)$ $\equiv$ State of the underlying Markov chain of the MAP.
\item $R_{s}(t)$ $\equiv$ Remaining service time.
\item $R_{v}(t)$ $\equiv$ Remaining vacation time.
\end{itemize}
$S(t)=0$ reflects the server's dormant status at time $t$.
Depending on vacation policy we have the following Markov process
\begin{equation*}
\begin{cases}
\{(q(t),S(t)),J(t)\}\cup\{\big(q(t),S(t),J(t),R_{s}(t)\big)\cup \big(q(t),K_{v}(t),J(t),R_{v}(t)\big)\}, &\mbox{for~ SV},\\
\{\big(q(t),S(t),J(t),R_{s}(t)\big)\cup \big(q(t),K_{v}(t),J(t),R_{v}(t)\big)\},&\mbox{for~ MV},
\end{cases}
\end{equation*}
with state space
\begin{equation*}
\begin{cases}
\{(n,0,i);0\leq n\leq h-1, 1\leq i\leq m\}\bigcup \{(n,r,i,u);n\geq 0, h\leq r\leq H,1\leq i \leq m,u\geq 0\}\bigcup\\ \{(n,k,i,u); 0\leq k\leq h-1, n\geq k,1\leq i \leq m, u\geq 0\}, &\mbox{for~ SV},\\
\{(n,r,i,u);n\geq 0, h\leq r\leq H,1\leq i \leq m,u\geq 0\}\bigcup\\ \{(n,k,i,u); 0\leq k\leq h-1, 1\leq i \leq m,n\geq k,u\geq 0\},&\mbox{for~ MV}.
\end{cases}
\end{equation*}

The state probabilities at time $t$ are defined as:
\begin{itemize}
\item $R_i{(n,0,t)} \equiv Pr{\{q(t)=n, S(t)=0, J(t)=i, u\leq R_{s}(t)\leq u+du\}}, 1\leq i\leq m , 0\leq n \leq h-1$, (for SV only).
\item $\xi_{i}{(n,r,u,t)}du \equiv Pr{\{q(t)=n,S(t)=r, J(t)=i,u \leq R_{s}(t) \leq u+du\}}, 1\leq i\leq m, n\geq 0~,~h\leq r\leq H$.
\item 	$\gamma_{i}(n,k,u,t)du\equiv~Pr{\{q(t)=n,K_{v}(t)=k, J(t)=i, u\leq R_{v}(t)\leq u+du \}}, 1\leq i\leq m, n\geq k~,~0\leq k\leq h-1$.
\end{itemize}
  In steady state, as $t\rightarrow\infty$ we have,

$R_i(n,0)=\displaystyle\ \lim_{t\to\infty}R_i{(n,0,t)},~~0\leq n\leq h-1, 1\leq i\leq m$, (exist only for SV),\\
 $\xi_{i}{(n,r,u)}=\displaystyle\ \lim_{t\to\infty}\xi_{i}{(n,r,u,t)},~~n\geq 0,~~ h\leq r\leq H,~1\leq i\leq m$,\\
 $\gamma_{i}{(n,k,u)}=\displaystyle\ \lim_{t\to\infty}\gamma_{i}{(n,k,u,t)},~~n\geq k,~~0\leq k\leq h-1,~1\leq i\leq m.$\\
 Further, we define
 \begin{itemize}
\item $R(n,0)$=($R_1(n,0),R_2(n,0),..,R_m(n,0)$),~~$0\leq n\leq h-1$.
\item $\xi(n,r,u)$=($\xi_1(n,r,u),\xi_2(n,r,u),...,\xi_m(n,r,u)$),~~$n\geq 0, h\leq r\leq H$.
\item $\gamma(n,k,u)$=$(\gamma_1(n,k,u),\gamma_2(n,k,u),...,\gamma_m(n,k,u))$,~~$n\geq k, 0\leq k\leq h-1$.
\end{itemize}
Following an analysis of the system at time $t$ and $t+dt$, the related steady state equations are obtained as follows:

\bea
0&=&(1-\epsilon)\bigg(R(0,0)C+\gamma
(0,0,0)\bigg),\label{eq9}\\
\nn0&=&(1-\epsilon)\bigg(R(n,0)C+R(n-1,0)D+\sum\limits_{k=0}^{n}\gamma(n,k,0)\bigg), 1\leq n\leq h-1, \label{eq10}\\
&&\\
\nn-\frac{d}{du}\xi{(0,h,u)}&=&\xi(0,h,u)C+(1-\epsilon)R(h-1,0)Ds_h{(u)}\\
 &&+\bigg(\sum\limits_{k=0}^{h-1}\gamma(a,k,0)+\sum\limits_{r=h}^{H}\xi(a,r,0)\bigg)s_h{(u)},\label{eq11}\\
-\frac{d}{du}\xi{(0,r,u)}&=&\xi(0,r,u)C+\bigg(\sum\limits_{k=0}^{h-1}\gamma(r,k,0)+\sum\limits_{r=h}^{H}\xi(r,j,0)\bigg)s_r{(u)}, h+1 \leq r \leq H, \label{eq12}\\
-\frac{d}{du}\xi{(n,r,u)}&=&\xi(n,r,u)C+\xi(n-1,r,u)D,~~ h\leq r \leq H-1,~~ n \geq 1,\label{eq13}\\
-\nn \frac{d}{du}\xi{(n,H,u)}&=&\xi{(n,H,u)}C+ \xi{(n-1,H,u)}D+\bigg(\sum\limits_{k=0}^{h-1}\gamma(n+H,k,0)\\
&&+\sum\limits_{r=h}^{H}\xi(n+H,r,0)\bigg)s_H{(u)},~n \geq 1,\label{eq14} \\
\nn-\frac{d}{du}\gamma(k,k,u)&=&\gamma(k,k,u)C+\bigg(\sum\limits_{r=h}^{H}\xi(k,r,0)+\epsilon\sum\limits_{j=0}^{k}\gamma(k,j,0)\bigg)v_k{(u)}, 0\leq k \leq h-1 \\ &&\label{eq15}\\
-\frac{d}{du}\gamma(n,k,u)&=&\gamma(n,k,u)C+\gamma(n-1,k,u)D,~~n\geq k+1,~~ 0\leq k \leq h-1.\label{eq16}
\eea
Further, for Re$(\theta)\geq 0$, we define,
\bea
\tilde{S}_r{(\theta)}=\int_{0}^{\infty}e^{-\theta u}dS_r{(u)}=\int_{0}^{\infty}e^{-\theta u}s_r{(u)}du,~~ h\leq r \leq H,\label{eq17}\\
\tilde \xi{(n,r,\theta)}=\int_{0}^{\infty}e^{-\theta u}\xi{(n,r,u)}du,~~h\leq r\leq H,~n\geq 0, \label{eq18}\\
\xi{(n,r)}\equiv \tilde \xi{(n,r,0)}=\int_{0}^{\infty}\xi(n,r,u)du,~~h\leq r\leq H,~ n\geq 0, \label{19}\\
\tilde{V}_k{(\theta)}=\int_{0}^{\infty}e^{-\theta u}dV_k{(u)}=\int_{0}^{\infty}e^{-\theta u}v_k{(u)}du,~~ 0 \leq k\leq h-1,\label{eq20}\\
\tilde \gamma(n,k,\theta)=\int_{0}^{\infty}e^{-\theta u}\gamma(n,k,u)du,~~0\leq k\leq h-1,~n\geq k, \label{eq21}\\
\gamma(n,k)\equiv \tilde \gamma(n,k,0)=\int_{0}^{\infty}\gamma(n,k,u)du,~~0\leq k\leq h-1,~~ n\geq k. \label{eq22}
\eea
 \begin{itemize}
\item $R(n,0)$=($R_1(n,0),R_2(n,0),...,R_m(n,0)$),~~$0\leq n\leq h-1$.
\item $\xi(n,r)$=($\xi_1(n,r),\xi_2(n,r),...,\xi_m(n,r)$),~~$n\geq 0, h\leq r\leq H$.
\item $\gamma(n,k)$=$(\gamma_1(n,k),\gamma_2(n,k),...,\gamma_m(n,k))$,~~$n\geq k, 0\leq k\leq h-1$.
\end{itemize}
Here the probability ($R_i(n,0)$)~\{$\xi_i{(n,r)}$\}~[$\gamma_i(n,k)$] denotes that (queue size is $n$ and the sever is dormant, and the arrival process is in phase $i$, $0\leq n\leq h-1$) \{queue size is $n$ and $r$ customers are being serviced, and the arrival process is in phase $i$, $h\leq r\leq H$, $n\geq 0$\} [queue size is $n$ and $k-th$ vacation type performs by the server, and the arrival process is in phase $i$, $0\leq k\leq h-1$, $n\geq k$] at arbitrary epoch.\\
If we multiply $e^{-\theta u}$ in (\ref{eq11})-(\ref{eq16}) and integrate with respect to variable $u$ with limits 0 to $\infty$, we get
\bea
\nn-\theta\tilde\xi{(0,h,\theta)}&=& \tilde\xi(0,h,\theta)C+(1-\epsilon) R{(h-1,0)}D\tilde{S_h}(\theta)+\bigg(\sum\limits_{k=0}^{h-1}\gamma(h,k,0)+\sum\limits_{r=h}^{H}\xi(h,r,0)\bigg)\tilde{S_h}(\theta)\\
&&-\xi{(0,h,0)},\label{eq23}\\
\nn-\theta\tilde\xi(0,r,\theta)&=&\tilde\xi(0,r,\theta)C+\bigg(\sum\limits_{k=0}^{h-1}\gamma(r,k,0)+\sum\limits_{r=h}^{H}\xi(r,j,0)\bigg)\tilde{S_r}(\theta)\\
&&-\xi(0,r,0),~~ h+1\leq r\leq H, \label{eq24}\\
-\theta\tilde\xi(n,r,\theta)&=&\tilde\xi(n,r,\theta)C+\tilde\xi(n-1,r,\theta)D-\xi(n,r,0),~~n\geq 1,h\leq r\leq H-1, \label{eq25}\\
-\theta\tilde\xi(n,H,\theta)&=&\tilde\xi(n,H,\theta)C+ \nn\tilde\xi(n-1,H,\theta)D+\bigg(\sum\limits_{k=0}^{h-1}\gamma(n+H,k,0)+\sum\limits_{r=h}^{H}\xi{(n+H,r,0)}\bigg)\tilde{S_H}(\theta)\\
&&-\xi(n,H,0),~n\geq 1,\label{eq26}\\
\nn-\theta\tilde{\gamma}(k,k,\theta)&=&\tilde{\gamma}(k,k,\theta)C+\bigg(\sum\limits_{r=h}^{H}\xi(k,r,0)+\epsilon\sum\limits_{j=0}^{k}\gamma(k,j,0)\bigg)\tilde{V}_k(\theta)-\gamma(k,k,0),~~0\leq k\leq h-1, \label{eq27}\\
&&\\
-\theta\tilde{\gamma}(n,k,\theta)&=&\tilde{\gamma}(n,k,\theta)C+\tilde \gamma(n-1,k,\theta)D-\gamma(n,k,0)~~n\geq k+1,~~0\leq k\leq h-1. \label{eq28}
\eea
Our aim is to perceive the probability vector of the joint probabilities of the queue content and server content (queue content and vacation type) at any time. However, direct analysis of these is quite challanging.The arbitrary epoch probabilities determine in terms of service (vacation) completion epoch probabilities after characterizing the system's state at the service (vacation) completion epoch. The following probabilities are defined in this context at the service (vacation) completion epoch while the arrival process is in phase $i$:

\bea
\nn \xi_i^{+}{(n,r)}&=&Pr\{\mbox{At the service completion epoch of a batch of size $r$, }\\
  &&\mbox{~~ queue size is $n$.}\},~~n\geq 0,~~h\leq r\leq H, \label{eq29}\\
\nn \xi_{i}^{+}(n)&=& Pr\{\mbox{At the service completion epoch, queue size is $n$}\}\\
&=&~~\sum\limits_{r=h}^{H}\xi_i^{+}{(n,r)},~~n\geq 0,~ \label{eq30}\\
\nn \gamma_i^{+}(n,k)&=&Pr\{\mbox{At $k-th$ type of vacation termination epoch,}\\
&& \mbox{queue size is $n$}\},~~0\leq k\leq h-1,~~n\geq k, \label{eq32}\\
\nn \gamma_i^{+}(n)&=&Pr\{\mbox{At vacation termination epoch, {queue size is $n$}}\\
&=&\sum\limits_{k=0}^{min(n,h-1)}\gamma_i^{+}(n,k),~~n\geq 0. \label{eq33}
\eea
Consequently, we obtain the probability vector.\\
$\xi^{+}(n,r)=(\xi_1^{+}(n,r),\xi_2^{+}(n,r),...,\xi_m^{+}(n,r)$),~~$n\geq 0,~~ h\leq r\leq H$,\\
$\gamma^{+}(n,k)=(\gamma_1^{+}(n,k),\gamma_2^{+}(n,k),...,\gamma_m^{+}(n,k))$,~~$n\geq k,~~ 0\leq k\leq h-1$,\\
$\xi^{+}(n)=(\xi_1^{+}(n),\xi_2^{+}(n),...,\xi_m^{+}(n))$,~~$n\geq 0$,\\
$\gamma^{+}(n)=(\gamma_1^{+}(n),\gamma_2^{+}(n),...,\gamma_m^{+}(n))$,~~$n\geq 0$.

\begin{lemma}\label{lma1}
The probability vectors $\xi^{+}{(n,r)}$, $\gamma^{+}(n,k)$, $\xi(n,r,0)$ and $\gamma(n,k,0)$ $(h\leq r\leq H,~0\leq k\leq h-1)$ are given by,
\bea
\xi^{+}{(n,r)}&=&\sigma{\xi(n,r,0)},\label{eq44}~~n\geq 0,\\
\gamma^{+}(n,k)&=&\sigma{\gamma(n,k,0)},\label{eq45}~~n\geq k,
\eea
where~ $\sigma^{-1}=\sum\limits_{m=0}^{\infty}\sum\limits_{r=h}^{H}\xi(m,r,0)$\bf e$+\sum\limits_{m=0}^{\infty}\sum\limits_{k=0}^{min(m,h-1)}\gamma(m,k,0)$\bf e$ $.
\end{lemma}
{\textit{Proof:}}
As a result of the fact that $\xi^{+}(n,r)$ and $\gamma^{+}(n,k)$ are proportional to $\xi(n,r,0)$ and $\gamma(n,k,0)$, respectively, applying Bayes' theorem and
$\sum\limits_{n=0}^{\infty}\sum\limits_{r=h}^{H}\xi^{+}(n,r)\mathbf{e} +\sum\limits_{n=0}^{\infty}\sum\limits_{k=0}^{min(n,h-1)}\gamma^{+}(n,k)\mathbf{e}$ =1 we get the desired outcome.

\begin{lemma}\label{lma2.1}
$R(n,0)D\mathbf{e}$=$\sum\limits_{m=0}^{n}\sum\limits_{k=0}^{m}\gamma(m,k,0)\mathbf{e}$
\end{lemma}
{\textit{Proof:}}
using (\ref{eq9}) and (\ref{eq10}) after some simplification we accomplish the intented outcome.
\begin{lemma}\label{lma2}
The expression for $\sigma^{-1}$ is
\bea
\sigma^{-1}=\frac{1-(1-\epsilon)\sum\limits_{n=0}^{h-1}R(n,0)\mathbf{e}}{\begin{split}s_H\sum\limits_{n=H+1}^{\infty}\big(\xi^{+}(n)+\gamma^{+}(n)\big)\mathbf{e}+\sum\limits_{n=h}^{H}\big(\xi^{+}(n)+\gamma^{+}(n)\big)\mathbf{e}s_n
\\+\sum\limits_{n=0}^{h-1}\big(\xi^{+}(n)x_n+(1-\epsilon)\gamma^{+}(n)s_h+\epsilon \gamma^{+}(n)x_n \big)\mathbf{e}\end{split}}.\label{eq46}
\eea
\end{lemma}
{\textit{Proof:}}
Post multiplying (\ref{eq23})-(\ref{eq28}) by column vector $\bf{e}$ and summing them using Lemma \ref{lma2.1} after some simplification we get
\bea
\nn\sum\limits_{m=0}^{\infty}\sum\limits_{r=h}^{H}\tilde\xi(m,r,\theta)+\sum\limits_{m=0}^{\infty}\sum\limits_{k=0}^{min(m,h-1)}\tilde \gamma(m,k,\theta)&=&\frac{1-\tilde{S}_H{(\theta)}}{\theta}\sum\limits_{n=H+1}^{\infty}\bigg(\sum\limits_{r=h}^{H}\xi(n,r,0)+\sum\limits_{k=0}^{h-1}\gamma(n,k,0)\bigg)\\
\nn&&+\sum\limits_{n=h}^{H}\bigg(\sum\limits_{r=h}^{H}\xi(n,r,0)+\sum\limits_{k=0}^{h-1}\gamma(n,k,0)\bigg)\frac{1-\tilde{S}_n{(\theta)}}{\theta} \label{eq48}\\
\nn&&+\sum\limits_{n=0}^{h-1}\bigg(\sum\limits_{r=h}^{H}\xi(n,r,0)+\epsilon\sum\limits_{k=0}^{n}\gamma(n,k,0)\bigg)\frac{1-\tilde{V}_n{(\theta)}}{\theta}\\
&&+(1-\epsilon)\frac{1-\tilde{S}_h{(\theta)}}{\theta}\sum\limits_{n=0}^{h-1}\sum\limits_{k=0}^{n}\gamma(n,k,0).
\eea
applying $\theta \rightarrow 0$ in (\ref{eq48}) and L'H\^{o}spital's rule, and $(1-\epsilon)\sum\limits_{n=0}^{h-1}R(n,0)\mathbf{e}+\sum\limits_{n=0}^{\infty}\sum\limits_{r=h}^{H}\xi(n,r)\mathbf{e}+\sum\limits_{n=0}^{\infty}\sum\limits_{k=0}^{min(n,h-1)}\gamma(n,k)\mathbf{e}$=1, after few simplification we get desired outcome.

Further, define a few necessary generating functions, which are as follows:
\bea
\tilde\Pi_i(z,y,\theta)&=&\sum\limits_{n=0}^{\infty}\sum\limits_{r=h}^{H}\tilde\xi(n,r,\theta)z^ny^r, \label{eq35}\\
\Pi_i^{+}(z,y)&=&\sum\limits_{n=0}^{\infty}\sum\limits_{r=h}^{H}\xi_i^{+}(n,r)z^ny^r, \label{eq36}\\
\Psi_i^{+}(z)&=&\sum\limits_{n=0}^{\infty}\sum\limits_{r=h}^{H}\xi_i^{+}(n,r)z^n=\sum\limits_{n=0}^{\infty}\xi_i^{+}(n)z^n=\Pi_i^{+}(z,1),\label{eq37}\\
\tilde O_i(z,y,\theta)&=&\sum\limits_{k=0}^{h-1}\sum\limits_{n=k}^{\infty}\tilde \gamma_i(n,k,\theta)z^ny^k, \label{eq38}\\
O_{i}^{+}(z,y)&=&\sum\limits_{k=0}^{h-1}\sum\limits_{n=k}^{\infty}\gamma_i^{+}(n,k)z^ny^k, \label{eq39}\\
\nn O_{i}^{+}(z)&=&\sum\limits_{k=0}^{h-1}\sum\limits_{n=k}^{\infty}\gamma_i^{+}(n,k)z^n=\sum\limits_{n=0}^{\infty}\sum\limits_{k=0}^{min(n,h-1)}\gamma_i^{+}(n,k)z^n=\sum\limits_{n=0}^{\infty}\gamma_i^{+}(n)z^n,\\
&& \label{eq40}
\eea
where, $|z|\leq 1$ and $|y|\leq 1$.
Consequently we have,\\
$\tilde\Pi(z,y,\theta)=(\tilde\Pi_1(z,y,\theta),\tilde\Pi_2(z,y,\theta),...,\tilde\Pi_m(z,y,\theta))$\\
$\Pi^{+}(z,y)=(\Pi_1^{+}(z,y),\Pi_2^{+}(z,y),...,\Pi_m^{+}(z,y))$\\
$\Psi^{+}(z)=(\Psi_1^{+}(z),\Psi_2^{+}(z),...,\Psi_m^{+}(z))$\\
$\tilde O(z,y,\theta)=(\tilde O_1(z,y,\theta),\tilde O_2(z,y,\theta),...,\tilde O_m(z,y,\theta)$\\
$O^{+}(z,y)=(O_1^{+}(z,y),O_2^{+}(z,y),...,O_m^{+}(z,y))$\\
$O^{+}(z)=(O_1^{+}(z),O_2^{+}(z),...,O_m^{+}(z))$
\begin{lemma}\label{lma4}
\bea
O^{+}(z)=\sum\limits_{n=0}^{\infty}\gamma^{+}(n)z^{n}=\sum\limits_{k=0}^{h-1}(\xi^{+}(k)+\epsilon \gamma^{+}(k))B^{(k)}(z)z^{k}\label{eq52}
\eea
\end{lemma}
{\textit{Proof:}}
Equations (\ref{eq27}) and (\ref{eq28}) are multiplied by the appropriate powers of $z$ and $y$ and added throughout the range of $n$ and $k$ we get
\bea
\nn \tilde O(z,y,\theta)(-\theta I-(C+Dz))&=&\sum\limits_{k=0}^{h-1}\bigg(\sum\limits_{r=h}^{H}\xi(k,r,0)+\epsilon \sum\limits_{j=0}^{k}\gamma(j,k,0)\bigg)\tilde{V}_k{(\theta)}z^{k}y^k\\
&&-\sum\limits_{k=0}^{h-1}\sum\limits_{n=k}^{\infty}\gamma(n,k,0)z^{n}y^k. \label{eq50}
\eea
If the eigenvalues of $-(C+Dz)$ are $\alpha_1(z),\alpha_2(z),..,\alpha_m(z)$ and $\delta_1(z),\delta_2(z),...,\delta_m(z)$ be the corresponding eigenvectors, then
\bea
-(C+Dz)\delta_i(z)=\alpha_i(z)\delta_i(z), 1\leq i\leq m.\label{eq50.1}
\eea
If we substitute $\theta=\alpha_i(z)$ in (\ref{eq50}) and post multiplying by $\delta_i(z)$, using Lemma \ref{lma1} and Lemma \ref{lma2} we get
\bea
\sum\limits_{k=0}^{h-1}\sum\limits_{n=k}^{\infty}\gamma^{+}(n,k)z^{n}y^k\delta_i(z)=\sum\limits_{k=0}^{h-1}(\xi^{+}(k)+\epsilon \gamma^{+}(k))\tilde V_k(\alpha_i(z))\delta_i(z)z^{k}y^k,\label{eq51}
\eea
Equation (\ref{eq51}) is true for all $\alpha_i(z)$, $1\leq i\leq m$, hence
\bea
\nn \sum\limits_{k=0}^{h-1}\sum\limits_{n=k}^{\infty}\gamma^{+}(n,k)z^{n}y^k=\sum\limits_{k=0}^{h-1}(\xi^{+}(k)+\epsilon \gamma^{+}(k))\Delta(z)diag\{\tilde V_{k}(\alpha_i(z))\}^m_{i=1}(\Delta(z))^{-1}z^{k}y^k.\\
&&\label{eq52.1}
\eea
where $\Delta(z)=(\delta_1(z),\delta_2(z),...,\delta_m(z))$ and $diag\{\tilde V_{k}(\alpha_i(z))\}^m_{i=1}$ is a diagonal matrix whose ($i,i$) entry is $\tilde V_{k}(\alpha_i(z))$, $i=1,2,...,m$.
\bea
\nn(B^{(k)}_l(x))_{i,j}&=&Pr\{\mbox{Given a departure at time 0 which left $k$ costumer in the queue,}\\
\nn &&\mbox{$k-th$ vacation type begins and the arrival process is in phase $i$}\\
\nn &&\mbox{at the end of the $k-th$ vacation type occurs no later than time $x$, with the arrival process is in}\\
\nn &&\mbox {phase $j$, and during the $k-th$ vacation type $l$ customers arrive}\}, 0\leq k\leq h-1.
\eea
Let $B^{(k)}(z)$ be the probability generating function of $B^{(k)}_l=(B^{(k)}_l(x))_{i,j}$, and hence,\\
$B^{(k)}(z)=\sum\limits_{l=0}^{\infty}B^{(k)}_lz^l=\int_{0}^{\infty}e^{-(C+Dz)t}v_{k}(t)dt=\Delta(z)diag\{\tilde V_{k}(\alpha_i(z))\}^m_{i=1}(\Delta(z))^{-1}$, $0\leq k\leq h-1$.\\
Hence, equation (\ref{eq52.1}) expresses as,
\bea
\sum\limits_{k=0}^{h-1}\sum\limits_{n=k}^{\infty}\gamma^{+}(n,k)z^{n}y^k=\sum\limits_{k=0}^{h-1}(\xi^{+}(k)+\epsilon \gamma^{+}(k))B^{(k)}(z)z^{k}y^k.\label{eq530}
\eea
substituting $y=1$ in (\ref{eq530}) we get desired result (\ref{eq52}).

When we multiply (\ref{eq23})-(\ref{eq26}) by the appropriate powers of $z$ and $y$ and add the results over the range of $n$ and $r$, we get

\bea
\nn\tilde\Pi(z,y,\theta)(-\theta I-(C+Dz))&=&(1-\epsilon)\sum\limits_{n=0}^{h-1}\sum\limits_{k=0}^{n}\gamma(n,k,0)\tilde D^{h-n}\tilde{S}_h{(\theta)}y^h\\
&&+\sum\limits_{r=h}^{H}\bigg(\sum\limits_{k=0}^{h-1}\gamma(r,k,0)+\sum\limits_{r=h}^{H}\xi(r,j,0)\bigg)\tilde{S}_r{(\theta)}y^r \label{eq53}\\
\nn&&+\sum\limits_{n=H+1}^{\infty}\bigg(\sum\limits_{k=0}^{h-1}\gamma(n,k,0)+\sum\limits_{r=h}^{H}\xi(n,r,0)\bigg)\tilde{S}_H{(\theta)}z^{n-H}y^H\\
\nn &&-\sum\limits_{n=0}^{\infty}\sum\limits_{r=h}^{H}\xi(n,r,0)z^{n}y^r,
\eea

where $\tilde{D}=(-C)^{-1}D$. Now we substitute $\theta=\alpha_i(z)$ in (\ref{eq53}) and post multiplying by $\delta_i(z)$, using Lemma \ref{lma1} and Lemma \ref{lma2} we get
\bea
\nn\sum\limits_{n=0}^{\infty}\sum\limits_{r=h}^{H}\xi^{+}(n,r)z^{n}y^r\delta_i(z)&=&(1-\epsilon)\sum\limits_{n=0}^{h-1}\gamma^{+}(n)\tilde D^{h-n}\tilde{S}_h{(\alpha_i(z))}\delta_i(z)y^h\\
\nn &&+\sum\limits_{r=h}^{H}\bigg(\gamma^{+}(r)+\xi^{+}(r)\bigg)\tilde{S}_r{(\alpha_i(z))}\delta_i(z)y^r\\
\nn &&+\sum\limits_{n=H+1}^{\infty}\bigg(\gamma^{+}(n)+\xi^{+}(n)\bigg)\tilde{S}_H{(\alpha_i(z))}\delta_i(z)z^{n-H}y^H.\\
&& \label{eq54}
\eea
Equation (\ref{eq54}) is true for all $\alpha_i(z)$, $1\leq i\leq m$, hence
\bea
\nn\sum\limits_{n=0}^{\infty}\sum\limits_{r=h}^{H}\xi^{+}(n,r)z^{n}y^r&=&(1-\epsilon)\sum\limits_{n=0}^{h-1}\gamma^{+}(n)\tilde D^{h-n}\Delta(z)diag\{\tilde{S}_h{(\alpha_i(z))}\}^{m}_{i=1}{(\Delta(z))}^{-1}y^h\\
&&+\sum\limits_{r=h}^{H}\bigg(\gamma^{+}(r)+\xi^{+}(r)\bigg)\Delta(z)diag\{\tilde{S}_r{(\alpha_i(z))}\}^{m}_{i=1}{(\Delta(z))}^{-1}y^r\\
\nn&&+\sum\limits_{n=H+1}^{\infty}\bigg(\gamma^{+}(n)+\xi^{+}(n)\bigg)\Delta(z)diag\{\tilde{S}_H{(\alpha_i(z))}\}^{m}_{i=1}{(\Delta(z))}^{-1}z^{n-H}y^H. \label{eq55}
\eea
Further, we define
\bea
\nn(A^{(r)}_l(x))_{i,j}&=&Pr\{\mbox{Given a departure at time 0 which left $r$ ($h\leq r\leq H$) costumer in the queue,}\\
\nn &&\mbox{and the arrival process is in phase $i$, next departure occurs no later than}\\
\nn &&\mbox{ time $x$ with the arrival process is in phase $j$, and during the}\\
\nn &&\mbox { service time of $r$ customers $l$ customers arrive}\}.
\eea

Let $A^{(r)}(z)$ be the probability generating function of $A^{(r)}_l=(A^{(r)}_l(x))_{i,j}$, and hence,\\
$A^{(r)}(z)=\sum\limits_{l=0}^{\infty}A^{(r)}_lz^l=\int_{0}^{\infty}e^{-(C+Dz)t}s_{r}(t)dt=\Delta(z)diag\{\tilde S_{r}(\alpha_i(z))\}^m_{i=1}(\Delta(z))^{-1}$, $h\leq r\leq H$.\\
Substituting $y=1$ in (\ref{eq54}) and using Lemma \ref{lma4} and equation (\ref{eq37}) we get the following result
\bea
\Psi^{+}(z)=\frac{\begin{split}\bigg\{\sum\limits_{n=0}^{h-1}\big[(\xi^{+}(n)+\epsilon \gamma^{+}(n))\big(B^{(n)}(z)-I\big)A^{(H)}(z)z^n\\+(1-\epsilon)\gamma^{+}(n)\big(\tilde D^{h-n}A^{(h)}(z)z^H-A^{(H)}(z)z^n\big)\big]\\+\sum\limits_{n=h}^{H-1}(\gamma^{+}(n)+\xi^{+}(n))\big(A^{(n)}(z)z^H -A^{(H)}(z)z^n\big)\bigg\}\end{split}}{z^HI-A^{(H)}(z)}.\label{eq55}
\eea
Now using (\ref{eq55}) in (\ref{eq54}) after algebraic manipulation the following expression is obtained,
\bea
\Pi^{+}(z,y)=\frac{\begin{split}\sum\limits_{n=0}^{h-1}\big[(1-\epsilon)\gamma^{+}(n)\big(\tilde D^{h-n}z^HA^{(h)}(z)y^{h}-A^{(H)}(z)z^ny^{H}\big)\\+(1-\epsilon)\gamma^{+}(n)\tilde D^{h-n}A^{(h)}(z)A^{(H)}(z)\big(y^{H}-y^{h}\big)\\
+y^{H}(\xi^{+}(n)+\epsilon \gamma^{+}(n))\big(B^{(n)}(z)-I\big)A^{(H)}(z)z^n\big]\\ +\sum\limits_{n=h}^{H-1}\big(\gamma^{+}(n)+\xi^{+}(n)\big)\big(z^Hy^{n}A^{(n)}(z)+\big(y^{H}-y^{n}\big)A^{(n)}(z)A^{(H)}(z)-y^{H}A^{(H)}(z)z^n\big)\end{split}}{z^HI-A^{(H)}(z)}.\label{eq56}
\eea
The above bivariate vector generating function given in (\ref{eq56}) contains $H$ unknowns vectors $\{\xi^{+}(n)\}^{H-1}_{n=0}$ i.e., total $mH$ unknowns $\{\xi^{+}_i(n)\}^{H-1}_{n=0}$, $1\leq i\leq m$ which has to be determine first.
  From (\ref{eq56}), the bivariate generating function $\Pi^{+}(z,y)$ has been represented in compact form, excluding the $H$ unknowns $\{\xi^{+}(n)\}_{n=0}^{H-1}$. Additionally, if we obtain $\xi^{+}(k)$ $(0\leq k\leq h-1)$ then from Lemma \ref{lma4} the probability vectors $\gamma^{+}(n,k)$ $(0\leq k\leq h-1)$ are known. As a result, in order to determine all of the probability vectors at the service (vacation) completion epoch, it is necessary to identify the unknowns $\{\xi^{+}(n)\}_{n=0}^{H-1}$.

 \subsection{Procedure of obtaining the unknowns $\xi^{+}(n)$ $(0\leq n\leq H-1)$} Let us take $\tilde S_r(\theta)$ and $\tilde V_k(\theta)$ both are rational function. Then each element of $A^{(r)}(z)$ ($h\leq r\leq H$) and $B^{(k)}(z)$ ($0\leq k\leq h-1$) are rational functions having same denominator say $d^{(r)}(z)$ ($0\leq r\leq H$) and $d^{(k)}(z)$ ($0\leq k\leq h-1$), respectively.
 Assign the ($i,j$)-th element of $A^{(r)}(z)$ say $\frac{f^{(r)}_{i,j}(z)}{d^{(r)}(z)}$, $1\leq i,j\leq m$, and the ($i,j$)-th element of $B^{(k)}(z)$ say $\frac{f^{(k)}_{i,j}(z)}{d^{(k)}(z)}$, $1\leq i,j\leq m$. Consequently, the ($i,j$)-th element of $z^HI-A^{(H)}(z)$ is
 \bea
 (z^HI-A^{(H)}(z))_{i,j}=\frac{\nu_{i,j}(z)}{d^{(H)}(z)}
 \eea
 where
\begin{equation*}
 \nu_{i,j}(z)=
 \begin{cases}
 z^H d^{(H)}(z)-f^{(H)}_{i,j}(z), i=j,\\
 -f^{(H)}_{i,j}(z), i\neq j.
 \end{cases}
 \end{equation*}
 Hence from (\ref{eq55}) we have $m$ system of equations in the matrix form
 \bea
 \Psi^{+}(z)M(z)=\Omega(z) \label{eq56.1}
 \eea
 where $(i,j)$-th entry of the matrix $M(z)$ is $\nu_{i,j}(z)$, and $\Omega(z)=(\Omega_1(z),\Omega_2(z),...,\Omega_m(z))^T$ is an $m\times1$ column matrix such that\\
 \bea
 \tiny
\nn \Omega_j(z)&=&\frac{\begin{split}{(1-\epsilon)\big(\prod_{r=0}^{h-1}d^{(r)}(z)\big)\big(\prod_{r=h+1}^{H-1}d^{(r)}(z)\big) \sum\limits_{n=0}^{h-1}\sum\limits_{i=1}^{m}\gamma_i^{+}(n)\bigg(z^H\sum\limits_{l=1}^{m}\tilde d_{i,l}^{(h-n)}f^{(h)}_{l,j}(z)d^{(H)}(z)}\\
 -z^nd^{(h)}(z)f^{(H)}_{i,j}(z)\bigg)+\prod_{r=h}^{H-1}d^{(r)}(z)\sum\limits_{n=0}^{h-1}\sum\limits_{i=1}^{m}(\xi^{+}_i(n)\\+\epsilon \gamma^{+}_i(n))\sum\limits_{l=1}^{m}u^{(n)}_{i,l}(z)f^{(H)}_{l,j}(z)z^n\prod^{h-1}_{r=0,r\neq n}d^{(r)}(z)\\
 +\sum\limits_{r=0}^{h-1}d^{(r)}(z)\sum\limits_{n=h}^{H-1}\sum\limits_{i=1}^{m}(\xi^{+}_i(n)+\gamma^{+}_i(n))\bigg(z^H d^{(H)}(z)f^{(n)}_{i,j}(z)\\-z^n f_{i,j}^{(H)}(z)d^{(n)}(z)\bigg)\prod^{H-1}_{r=h,r\neq n}d^{(r)}(z)
 \end{split}}{\prod^{H-1}_{r=0}d^{(r)}(z)},1\leq j\leq m.\\
&&\label{eq56.2}
 \eea
 where $u^{(k)}_{i,j}(z)=f^{(k)}_{i,j}(z)-d^{(k)}(z)$, $i\neq j$, $0\leq k\leq h-1$ and $u^{(k)}_{i,j}(z)=f^{(k)}_{i,j}(z)$, $i=j$, $0\leq k\leq h-1$. $\tilde d_{k,l}^{(h-n)}$ is the $(k,l)$-th element of $\tilde D^{h-n}$.
To solve the system of equation given in (\ref{eq56.1}) we apply the cramer's rule. The result is as follows:
\bea
\Psi_j^{+}(z)=\frac{|M_j(z)|}{|M(z)|},~~ 1\leq j\leq m\label{eq56.3}
\eea
\begin{equation*}
 [M_{j}(z)]_{k,l}=
 \begin{cases}
 \Omega_k(z), j=l,\\
 \nu_{l,k}(z), j\neq l.
 \end{cases}
 \end{equation*}
Suppose that $|M(z)|$ is a non-zero polynomial in variable $z$ must posses a nonzero coefficient of the power of $z$.
 It is clear to observe that $|z^HI-A^{(H)}(z)|=\frac{|M(z)|}{{(d^{(H)}(z))^m}}$ has  precisely $mH$ zeros in \{$z: |z| \leq 1$\} say $p_1,p_2,...,p_l$ with multiplicity $q_1,q_2,...,q_l$, respectively, \big(where ($l\leq mH-1$) and $\sum\limits_{i=1}^{l}q_i=(mH-1)$\big) and $p_H=1$ is a simple zero. Since $\Psi_j^{+}(z)$ is analytic in $|z|\leq 1$ therefore, these roots are also the roots of the numerator of $\Psi_j^{+}(z)$. Hence taking one component of $\Psi^{+}(z)$ say
 $\Psi_{j}^{+}(z)$ ($1\leq j\leq m$), we are led to $mH-1$ equation
 \bea
 \bigg[\frac{d^{i-1}}{dz^{i-1}}|M_j(z)|\bigg]_{z=p_x}=0,~~1\leq x\leq l~~ \&~~ 1\leq i\leq q_j,\label{eq62a}
 \eea
 where $\frac{d^0}{dz^0}h(z)=h(z)$.\\
One more equation is obtained by the normalization condition $\Psi^{+}(1)\mathbf{e}+O^{+}(1)\mathbf{e}=1$, i.e,
\bea
\sum\limits_{j=1}^{m}\bigg[{\frac{d}{dz}|M_j(z)|}\bigg]_{z=1}+\bigg[{\frac{d}{dz}|M(z)|}\bigg]_{z=1}\sum\limits_{k=0}^{h-1}(\xi^{+}(k)+\epsilon \gamma^{+}(k))e=\bigg[{\frac{d}{dz}|M(z)|}\bigg]_{z=1}.\label{eq62.11new}
\eea
Solving (\ref{eq62a}) and (\ref{eq62.11new}) together we get $mH$ unknowns $\xi_j^{+}(n)$ ($1\leq j\leq m, 0\leq n\leq H-1$).
\begin{theorem}
The probability vectors of the joint probability of queue and server content are given by
\bea
\xi^{+}{(n,h)}&=& \bigg((1-\epsilon)\sum\limits_{m=0}^{h-1}\gamma^{+}(m)\tilde D^{h-m}+\gamma^{+}(h)+\xi^{+}{(h)}\bigg)A^{(h)}_n,\label{eq62}\\
\xi^{+}{(n,r)}&=& \bigg(\gamma^{+}(r)+\xi^{+}{(r)}\bigg)A^{(r)}_n,~~{h+1}\leq r\leq {H-1}.\label{eq63}
\eea
\end{theorem}
{\textit{Proof:}}
Using (\ref{eq36}) in (\ref{eq56}) collecting the multiples of $y^r$ $(h\leq r\leq H-1)$, we obtain
\bea
{\mbox{multiple of $y^h$}} :\sum\limits_{n=0}^{\infty}\xi_{n,h}^{+}z^n&=&\bigg((1-\epsilon)\sum\limits_{m=0}^{h-1}\gamma^{+}(m) \tilde D^{h-m}+\gamma^{+}(h)+\xi^{+}{(h)}\bigg)A^{(h)}(z)\label{eq66}\\
{\mbox{multiple of $y^r$}}:\sum\limits_{n=0}^{\infty}\xi_{n,r}^{+}z^n&=&\bigg(\gamma^{+}(r)+\xi^{+}{(r)}\bigg)A^{(r)}(z), h+1\leq r\leq H-1\label{eq67}
\eea
Accumulating  the coefficients of $z^n$,  from both side of (\ref{eq66}) and (\ref{eq67}), we readily obtain desired result (\ref{eq62}) and (\ref{eq63}).
Now our current objective is to collect the remaining probability vectors $\xi^{+}{(n,H)}$ $(n\geq 0)$. Towards this end, using (\ref{eq36}) in (\ref{eq56}) and then collecting the multiple of $y^H$ we get
\bea
\sum\limits_{n=0}^{\infty}\xi^{+}{(n,H)}z^n=\frac{\begin{split}A^{(H)}(z)\bigg\{\sum\limits_{n=0}^{h-1}\bigg[(\xi^{+}(n)+\epsilon \gamma^{+}(n))(B^{(n)}(z)-I)z^n\\+(1-\epsilon) \gamma^{+}(n)(\tilde D^{h-n}A^{(h)}(z)-z^n I)\bigg]\\+\sum\limits_{n=h}^{H-1}(\gamma^{+}(n)+\xi^{+}(n))(A^{(n)}(z)-z^n I)\bigg\}\end{split}}{z^H I-A^{(H)}(z)}.\label{eq68}
\eea
Assign a symbol $\sum\limits_{n=0}^{\infty}\xi^{+}{(n,H)}z^n$ as $\pounds^{+}(z)=(\pounds_{1}^{+}(z),\pounds_{2}^{+}(z),...,\pounds_{m}^{+}(z))$, and we replace $\Psi^{+}(z)$ and $\Omega_j(z)$ by $\pounds^{+}(z)$ and $\Theta_j(z)$, respectively, where $\Theta_j(z)$, ($1\leq j\leq m$) is given by
\bea
\nn \Theta_j(z)&=&\frac{\begin{split}{(1-\epsilon)\big(\prod_{r=0}^{h-1}d^{(r)}(z)\big)\big(\prod_{r=h+1}^{H-1}d^{(r)}(z)\big)\sum\limits_{n=0}^{h-1}\sum\limits_{i=1}^{m}\gamma_i^{+}(n)\bigg(\sum\limits_{w=1}^{m}\sum\limits_{l=1}^{m}\tilde d_{i,l}^{(h-n)}f^{(h)}_{l,w}(z)f^{(H)}_{w,j}(z)}\\
 -z^nd^{(h)}(z)f^{(H)}_{i,j}(z)\bigg)+\prod_{r=h}^{H-1}d^{(r)}(z)\sum\limits_{n=0}^{h-1}\sum\limits_{i=1}^{m}(\xi^{+}_i(n)\\+\epsilon \gamma^{+}_i(n))\sum\limits_{l=1}^{m}u^{(n)}_{i,l}(z)f^{(H)}_{l,j}(z)z^n\prod^{h-1}_{r=0,r\neq n}d^{(r)}(z)\\
 +\prod_{r=0}^{h-1}d^{(r)}(z)\sum\limits_{n=h}^{H-1}\sum\limits_{i=1}^{m}(\xi^{+}_i(n)+\gamma^{+}_i(n))\sum\limits_{w=1}^{m}\bigg(f^{(n)}_{i,w}(z)f^{(H)}_{w,j}(z)\\-z^n f_{i,j}^{(H)}(z)d^{(n)}(z)\bigg)\prod^{H-1}_{r=h,r\neq n}d^{(r)}(z)
 \end{split}}{\prod^{H-1}_{r=0}d^{(r)}(z)}, 1\leq j \leq m.\\
&&\label{eq56.2a}
\eea
As a result, $\pounds^{+}_{j}(z)$ is expressed as
\bea
\pounds^{+}_{j}(z)=\frac{|N_j(z)|}{|M(z)|},~~ 1\leq j\leq m, \label{eq69}
\eea
where,
\begin{equation*}
 [N_{j}(z)]_{k,l}=
 \begin{cases}
 \Theta_k(z), j=l,\\
 \nu_{l,k}(z), j\neq l.
 \end{cases}
 \end{equation*}
It is assumed that degree of $|N_j(z)|$ is $\hat{d_j}$ and degree of $|M(z)|$ is $\hat{d}$.

The roots of $|N_j(z)|$ of modules with more than one must be known in order to derive the probability vectors $\xi^{+}{(n,H)}$ $(n\geq 0)$.

Since (\ref{eq68}) is analytic, in $|z|\leq 1$, the roots of $|M(z)|$  lying in $|z|\leq 1$ are also the roots of $|N_j(z)|$, hence,  the roots lying in $|z|\leq 1$ can not be used to calculate $\xi^{+}{(n,H)}$ $(n\geq 0)$.
let us assume $\beta_1,\beta_2,...,\beta_l$ are the zeros of $|M(z)|$ of modules greater than one having multiplicity $\eta_1,\eta_2,...,\eta_l$, respectively, and $\sum\limits_{j=1}^{l}\eta_j<\hat{d} $. Here, two cases arrise\\\\
\emph{Case A:} $\hat{d}\leq \hat{d_j}$\\
now we apply the partial fraction method on (\ref{eq69}), $\pounds^{+}_{j}(z)$ can be written as,\\
\bea
\pounds^{+}_{j}(z)=\sum\limits_{i=0}^{\hat{d_j}-\hat{d}}\varrho_{i,j}z^i+\sum\limits_{w=1}^{l}\sum\limits_{i=1}^{\eta_{w}}\frac{B_{i,w,j}}{(z-\beta_{w})^{\eta_{w}-i+1}},\label{eq69a}
\eea
where
\bea
\nn B_{i,w,j}=\frac{1}{(i-1)!}\bigg[\frac{d^{i-1}}{dz^{i-1}}\bigg(\frac{|N_j(z)|\frac{d^{\eta_w}}{dz^{\eta_w}}(z-\beta_w)^{\eta_w}}{\frac{d^{\eta_w}}{dz^{\eta_w}}|M(z)|}\bigg)\bigg]_{z=\beta_{w}}, w=1,2,...,l,~i=1,2,...,\eta_w,\\
\nn j=1,2,...,m.
\eea
Collecting the coefficients of $z^n$ $(n\geq 0)$ from both side of (\ref{eq69a}) for ($1\leq j\leq m$) we obtain,
\begin{equation*}
\xi^{+}_{j}{(n,H)}=
\begin{cases}
\varrho_{n,j}+\sum\limits_{w=1}^{l}\sum\limits_{i=1}^{\eta_w}\frac{B_{i,w,j}}{(-1)^{\eta_w-i+1}\beta_w^{\eta_{w}+n-i+1}}\binom{{\eta_{w}-i+n}}{{\eta_w-i}}, &\mbox{$0\leq n\leq \hat{d_j}-\hat{d}$},\\
\sum\limits_{w=1}^{l}\sum\limits_{i=1}^{\eta_w}\frac{B_{i,w,j}}{(-1)^{\eta_w-i+1}\beta_w^{\eta_{w}+n-i+1}}\binom{{\eta_{w}-i+n}}{{\eta_w-i}},&\mbox{$n>\hat{d_j}-\hat{d}$}.
\end{cases}
\end{equation*}
\emph{Case B}: $\hat{d}>\hat{d_j}$\\
We remove the first summation term on the right hand side of (\ref{eq69a}) then,
\bea
\xi^{+}_{j}{(n,H)}=\sum\limits_{w=1}^{l}\sum\limits_{i=1}^{\eta_w}\frac{B_{i,w,j}}{(-1)^{\eta_w-i+1}\beta_w^{\eta_{w}+n-i+1}}\binom{{\eta_{w}-i+n}}{{\eta_w-i}},~~n\geq 0.
\eea
\begin{theorem}
Arbitrary epoch probability vectors  are given by,
\bea
R(n,0)&=&\frac{\sum\limits_{m=0}^{n}\gamma^{+}(m)\tilde D^{n-m}(-C)^{-1}}{E},~~0\leq n\leq h-1~~ ({exist~only~for~ SV)}\label{eq73}\\
\nn \xi{(0,h)}&=& (1-\epsilon)R(h-1,0)D+\bigg(\frac{\xi^{+}(h)+\gamma^{+}(h)-\xi^{+}(0,h)}{E}\bigg)(-C)^{-1},~~ n\geq 0,\\
&& \label{eq74}\\
\xi{(0,r)}&=&\bigg(\frac{\xi^{+}(r)+\gamma^{+}(r)-\xi^{+}(0,r)}{E}\bigg)(-C)^{-1},~~ n\geq 0,~ h+1\leq r\leq H-1, \label{eq75}\\
\xi{(n,r)}&=&\bigg(\xi{(n-1,r)}D-\frac{\xi^{+}(n,r)}{E}\bigg)(-C)^{-1},~~ n\geq 1, \label{eq76}\\
\nn \xi{(n,H)}&=&\bigg(\xi{(n-1,H)}D+\frac{\xi^{+}(n+H)+\gamma^{+}(n+H)-\xi^{+}(n,H)}{E}\bigg)(-C)^{-1},~n\geq 0,\\
&& \label{eq76.a}\\
\gamma(k,k)&=&\bigg(\frac{\xi^{+}(k)+\epsilon \gamma^{+}(k)-\gamma^{+}(k,k)}{E}\bigg)(-C)^{-1}, ~0\leq k\leq h-1, \label{eq77}\\
\gamma(n,k)&=&\bigg(\gamma(n-1,k)-\frac{\gamma^{+}(n,k)}{E}\bigg)(-C)^{-1},~~ n\geq k+1, ~0\leq k\leq h-1. \label{eq77}
\eea
where $E=\hat{w}+(1-\epsilon)\sum\limits_{n=0}^{h-1}\sum\limits_{m=0}^{n}\gamma^{+}(m)\tilde D^{(n-m)}(-C)^{-1}\mathbf{e}$,\\ $\hat{w}=s_H\sum\limits_{n=H+1}^{\infty}\big(\xi^{+}(n)+\gamma^{+}(n)\big)\mathbf{e}+\sum\limits_{n=h}^{H}\big(\xi^{+}(n)+\gamma^{+}(n)\big)\mathbf{e} s_n+\sum\limits_{n=0}^{h-1}\big(\xi^{+}(n)\mathbf{e}x_n+(1-\epsilon) \gamma^{+}(n)\mathbf{e}s_h+\epsilon \gamma^{+}(n)\mathbf{e}x_n \big)$.
\end{theorem}
{\textit{Proof:}}
	Dividing $\sigma^{-1}$ in equation (\ref{eq9}) and (\ref{eq10}) after simple algebraic manipulation we obtain (\ref{eq73}). Futher, taking $\theta \rightarrow 0$ in (\ref{eq11})-(\ref{eq17}) and then diving $\sigma^{-1}$ after simple algebraic manipulation we get desired outcome (\ref{eq74})-(\ref{eq77}).
\section{Marginal Probabilities}\label{sec4} Some marginal probabilities are given as follows:
 \begin{enumerate}
 \item \emph{Queue length distribution} is given by
 \begin{equation*}
 P_n^{queue}=
 \begin{cases}
 (1-\epsilon) R(n,0)\mathbf{e}+\sum\limits_{r=h}^{H}\xi{(n,r)}\mathbf{e}+\sum\limits_{k=0}^{min(n,h-1)}\gamma(n,k)\mathbf{e},~0\leq n\leq h-1,\\
 \sum\limits_{r=h}^{H}\xi{(n,r)}\mathbf{e}+\sum\limits_{k=0}^{min(n,h-1)}\gamma(n,k)\mathbf{e},~n\geq h.
 \end{cases}
 \end{equation*}
 \item \emph{The probability that the server is in a dormant state} ($P^{dor}$)=$\sum\limits_{n=0}^{h-1}R{(n,0)}\mathbf{e}$.
 \item \emph{Probability of $r$ customers with the server} ($P_r^{ser}$)=$\sum\limits_{n=0}^{\infty}\xi{(n,r)}\mathbf{e}$,~$h\leq r\leq H$.
 \item \emph{Probability of server's vacation is $k-th$ vacation type} ($\gamma_{vac}^{[k]}$)=$\sum\limits_{n=k}^{\infty}\gamma(n,k)\mathbf{e}$,~$0\leq k\leq h-1$.
 \item \emph{The probability that the server is busy} ($P_{busy}$)=$\sum\limits_{r=h}^{H}\sum\limits_{n=0}^{\infty}\xi{(n,r)}\mathbf{e}$.
 \item \emph{The probability that the server is on vacation} ($\gamma_{vac}$)=$\sum\limits_{k=0}^{h-1}\sum\limits_{n=k}^{\infty}\gamma(n,k)\mathbf{e}$.
 \end{enumerate}
 \section{Performance measure}\label{sec5} Performance measure is collected for observing the system performance. It helps the system manager for observe the system behavior so that he can modify the system getting for more efficient result. In this section, we collect significant performance measures which are as follows.

 \begin{enumerate}
  \item \emph{The expected number in the queue} ($L_q$) = $(1-\epsilon)\sum\limits_{n=0}^{h-1}nR(n,0)\mathbf{e}+\sum\limits_{n=0}^{\infty}\sum\limits_{r=h}^{H}n\xi{(n,r)}\mathbf{e}+\sum\limits_{k=0}^{h-1}\sum\limits_{n=k}^{\infty}n\gamma(n,k)\mathbf{e}$=$(1-\epsilon)\sum\limits_{n=0}^{h-1}nP_n^{queue}+\sum\limits_{n=h-\epsilon h}^{\infty}nP_n^{queue}$.
  \item \emph{The expected number in the system}  ($L_s$) = $(1-\epsilon)\sum\limits_{n=0}^{h-1}nR(n,0)\mathbf{e}+\sum\limits_{n=0}^{\infty}\sum\limits_{r=h}^{H}(n+r)\xi{(n,r)}\mathbf{e}+\sum\limits_{k=0}^{h-1}\sum\limits_{n=k}^{\infty}n\gamma(n,k)\mathbf{e}$.
 \item \emph{The expected waiting time in the queue} ($W_q$) =$\frac{L_q}{\lambda}$.
 \item \emph{The expected waiting time in the system} ($W_s$) =$\frac{L_s}{\lambda}$.
 \item \emph{Expected number with the server when server is busy} ($L^{ser}$) =$\sum\limits_{r=h}^{H}(rP_r^{ser}/P_{busy})$.
 \item \emph{Expected vacation type when server is on vacation} ($L^{vac})=\sum\limits_{k=0}^{h-1}(k\gamma_{vac}^{[k]}/\gamma_{vac})$.
 \end{enumerate}
\section{Numerical results}\label{sec6} The main objective for presenting this section is to validate the mathematical results presented in the previous section with some numerical results. These results are displayed in the tabular and graphical form, considering the service (vacation) time distribution as phase (PH) - type, which is also known as the representation as ($\alpha, T$), where $\alpha$ is a row vector of order $1\times n$, and $T$ is a square matrix of order $n$. The joint probabilities with predefined notations are presented for $MAP/G_r^{(5,9)}
/1$ queue with queue size dependent SV (MV) in the tabular form (Table \ref{maptab1} - Table \ref{maptab8}). The service time of each batch and vacation time of the server follow the PH-type distribution. The other input parameters are given below. The $MAP$ is represented
by the matrices $C=$
$\begin{pmatrix}
  -91.8125 & 14.1250 \\
  49.4375 & -77.6875
\end{pmatrix}$
and $D=$
$\begin{pmatrix}
  49.4375 & 28.2500\\
  7.0625 & 21.1875
\end{pmatrix}$. The service time of each batch under service follow the Erlang ($E_3$) distribution having PH-type representation ($\alpha_r, T_r$), where $T_r$= $\begin{pmatrix}
  -\mu_r & \mu_r & 0.0\\
  0.0 & -\mu_r & \mu_r\\
  0.0 & 0.0 & -\mu_r
\end{pmatrix}$,  $\mu_r=\frac{r(5.2)}{2}$ , $\alpha_r$= $\begin{pmatrix}
  1.0 & 0.0 & 0.0, \end{pmatrix}$, $5\leq r\leq 9$. The vacation time of the server follows $E_2$ distribution, having PH-type representation ($\alpha_k, T_k$), $T_k$= $\begin{pmatrix}
  -\nu_k & \nu_k\\
  0.0 & -\nu_k
\end{pmatrix}$
, $\nu_k= (k+1.0)^{2}0.5$, $0\leq k\leq 4$.
  $\alpha_k$= $\begin{pmatrix}
 1.0 & 0.0
\end{pmatrix}$, $0\leq k \leq 4$. $\xi$=[0.57143, 0.42857],  $\lambda=56.50$.


\begin{table}[H]
\tiny
  \centering
  \caption{Joint probabilities (queue size and server content) at service completion epoch for SV}
  \smallskip
    \begin{tabular}{|c|c|c|c|c|c|c|c|c|c|c|}
    \hline
    {$n$} & {$\xi_{1}^{+}(n,5)$} & {$\xi_{2}^{+}(n,5)$} & {$\xi_{1}^{+}(n,6)$} & {$\xi_{2}^{+}(n,6)$} & {$\xi_{1}^{+}(n,7)$} & {$\xi_{2}^{+}(n,7)$} & {$\xi_{1}^{+}(n,8)$} & {$\xi_{2}^{+}(n,8)$} & {$\xi_{1}^{+}(n,9)$} & {$\xi_{2}^{+}(n,9)$} \\
    \hline
    0     & 0.00243 & 0.00173 & 0.00227 & 0.00162 & 0.00256 & 0.00183 & 0.00274 & 0.00196 & 0.00283 & 0.00203 \\
    1     & 0.00430 & 0.00317 & 0.00371 & 0.00275 & 0.00390 & 0.00289 & 0.00389 & 0.00289 & 0.00626 & 0.00459 \\
    2     & 0.00508 & 0.00379 & 0.00405 & 0.00304 & 0.00395 & 0.00297 & 0.00368 & 0.00278 & 0.00888 & 0.00659 \\
    3     & 0.00500 & 0.00376 & 0.00369 & 0.00278 & 0.00334 & 0.00253 & 0.00291 & 0.00221 & 0.01035 & 0.00773 \\
    4     & 0.00443 & 0.00334 & 0.00302 & 0.00229 & 0.00254 & 0.00193 & 0.00207 & 0.00158 & 0.01086 & 0.00814 \\
    5     & 0.00366 & 0.00277 & 0.00231 & 0.00175 & 0.00181 & 0.00138 & 0.00137 & 0.00105 & 0.01074 & 0.00807 \\
    6     & 0.00288 & 0.00219 & 0.00168 & 0.00128 & 0.00122 & 0.00093 & 0.00087 & 0.00066 & 0.01027 & 0.00772 \\
    7     & 0.00219 & 0.00166 & 0.00118 & 0.00090 & 0.00080 & 0.00061 & 0.00053 & 0.00040 & 0.00964 & 0.00725 \\
    8     & 0.00162 & 0.00123 & 0.00080 & 0.00061 & 0.00050 & 0.00039 & 0.00031 & 0.00024 & 0.00897 & 0.00675 \\
    9     & 0.00117 & 0.00089 & 0.00054 & 0.00041 & 0.00031 & 0.00024 & 0.00018 & 0.00014 & 0.00832 & 0.00626 \\
    10    & 0.00083 & 0.00063 & 0.00035 & 0.00027 & 0.00019 & 0.00015 & 0.00010 & 0.00008 & 0.00773 & 0.00581 \\
    31    & 0.00000 & 0.00000 & 0.00000 & 0.00000 & 0.00000 & 0.00000 & 0.00000 & 0.00000 & 0.00294 & 0.00220 \\
    32    & 0.00000 & 0.00000 & 0.00000 & 0.00000 & 0.00000 & 0.00000 & 0.00000 & 0.00000 & 0.00285 & 0.00214 \\
    33    & 0.00000 & 0.00000 & 0.00000 & 0.00000 & 0.00000 & 0.00000 & 0.00000 & 0.00000 & 0.00277 & 0.00208 \\
    151   & 0.00000 & 0.00000 & 0.00000 & 0.00000 & 0.00000 & 0.00000 & 0.00000 & 0.00000 & 0.00047 & 0.00035 \\
    152   & 0.00000 & 0.00000 & 0.00000 & 0.00000 & 0.00000 & 0.00000 & 0.00000 & 0.00000 & 0.00046 & 0.00035 \\
    153   & 0.00000 & 0.00000 & 0.00000 & 0.00000 & 0.00000 & 0.00000 & 0.00000 & 0.00000 & 0.00046 & 0.00034 \\
    301   & 0.00000 & 0.00000 & 0.00000 & 0.00000 & 0.00000 & 0.00000 & 0.00000 & 0.00000 & 0.00006 & 0.00004 \\
    302   & 0.00000 & 0.00000 & 0.00000 & 0.00000 & 0.00000 & 0.00000 & 0.00000 & 0.00000 & 0.00006 & 0.00004 \\
    303   & 0.00000 & 0.00000 & 0.00000 & 0.00000 & 0.00000 & 0.00000 & 0.00000 & 0.00000 & 0.00006 & 0.00004 \\
    372   & 0.00000 & 0.00000 & 0.00000 & 0.00000 & 0.00000 & 0.00000 & 0.00000 & 0.00000 & 0.00002 & 0.00001 \\
    373   & 0.00000 & 0.00000 & 0.00000 & 0.00000 & 0.00000 & 0.00000 & 0.00000 & 0.00000 & 0.00002 & 0.00001 \\
    442   & 0.00000 & 0.00000 & 0.00000 & 0.00000 & 0.00000 & 0.00000 & 0.00000 & 0.00000 & 0.00001 & 0.00001 \\
    443   & 0.00000 & 0.00000 & 0.00000 & 0.00000 & 0.00000 & 0.00000 & 0.00000 & 0.00000 & 0.00001 & 0.00000 \\
    444   & 0.00000 & 0.00000 & 0.00000 & 0.00000 & 0.00000 & 0.00000 & 0.00000 & 0.00000 & 0.00001 & 0.00000 \\
    459   & 0.00000 & 0.00000 & 0.00000 & 0.00000 & 0.00000 & 0.00000 & 0.00000 & 0.00000 & 0.00001 & 0.00000 \\
    460   & 0.00000 & 0.00000 & 0.00000 & 0.00000 & 0.00000 & 0.00000 & 0.00000 & 0.00000 & 0.00001 & 0.00000 \\
    461   & 0.00000 & 0.00000 & 0.00000 & 0.00000 & 0.00000 & 0.00000 & 0.00000 & 0.00000 & 0.00001 & 0.00000 \\
    $\geq$ 462   & 0.00000 & 0.00000 & 0.00000 & 0.00000 & 0.00000 & 0.00000 & 0.00000 & 0.00000 & 0.00000 & 0.00000 \\
    \hline
    Total & 0.0353	& 0.0265 &	0.0242 &	0.0181& 0.021 &	0.0160	& 0.0188 &	0.0141 &	0.363 &	0.272 \\
    \hline
   \end{tabular}%
  \label{maptab1}%
\end{table}%

\begin{table}[H]
\tiny
  \centering
  \caption{Joint probabilities (queue size and vacation type) at vacation termination epoch for SV}
  \smallskip
    \begin{tabular}{|p{0.8cm}|p{0.8cm}|p{0.8cm}|p{0.8cm}|p{0.8cm}|p{0.8cm}|p{0.8cm}|p{0.8cm}|p{0.8cm}|p{0.8cm}|p{0.8cm}|p{0.8cm}|}
    \hline
    {$n$} & {$\gamma_1^{+}(n,0)$} & {$\gamma_2^{+}(n,0)$} & {$\gamma_1^{+}(n,1)$} & {$\gamma_{2}^{+}(n,1)$} & {$\gamma_{1}^{+}(n,2)$} & {$\gamma_{2}^{+}(n,2)$} & {$\gamma_{1}^{+}(n,3)$} & {$\gamma_{2}^{+}(n,3)$} & {$\gamma_{1}^{+}(n,4)$} & {$\gamma_{2}^{+}(n,4)$} & {$P_{n}^{queue+}$} \\
    \hline
    0     & 0.00000 & 0.00000 & 0.00000 & 0.00000 & 0.00000 & 0.00000 & 0.00000 & 0.00000 & 0.00000 & 0.00000 & 0.02200 \\
    1     & 0.00001 & 0.00001 & 0.00010 & 0.00007 & 0.00000 & 0.00000 & 0.00000 & 0.00000 & 0.00000 & 0.00000 & 0.03853 \\
    2     & 0.00001 & 0.00001 & 0.00018 & 0.00013 & 0.00049 & 0.00035 & 0.00000 & 0.00000 & 0.00000 & 0.00000 & 0.04598 \\
    3     & 0.00001 & 0.00001 & 0.00025 & 0.00019 & 0.00084 & 0.00062 & 0.00124 & 0.00089 & 0.00000 & 0.00000 & 0.04834 \\
    4     & 0.00002 & 0.00001 & 0.00031 & 0.00023 & 0.00108 & 0.00080 & 0.00193 & 0.00143 & 0.00218 & 0.00156 & 0.04975 \\
    5     & 0.00002 & 0.00002 & 0.00036 & 0.00027 & 0.00124 & 0.00093 & 0.00225 & 0.00168 & 0.00301 & 0.00224 & 0.04692 \\
    6     & 0.00002 & 0.00002 & 0.00041 & 0.00031 & 0.00134 & 0.00100 & 0.00234 & 0.00175 & 0.00312 & 0.00234 & 0.04236 \\
    7     & 0.00003 & 0.00002 & 0.00045 & 0.00033 & 0.00139 & 0.00104 & 0.00228 & 0.00171 & 0.00288 & 0.00217 & 0.03745 \\
    8     & 0.00003 & 0.00002 & 0.00048 & 0.00036 & 0.00140 & 0.00105 & 0.00213 & 0.00160 & 0.00250 & 0.00188 & 0.03286 \\
    9     & 0.00003 & 0.00002 & 0.00050 & 0.00037 & 0.00138 & 0.00103 & 0.00193 & 0.00146 & 0.00207 & 0.00157 & 0.02883 \\
    10    & 0.00004 & 0.00003 & 0.00052 & 0.00039 & 0.00133 & 0.00100 & 0.00172 & 0.00130 & 0.00168 & 0.00127 & 0.02540 \\
    31    & 0.00007 & 0.00005 & 0.00038 & 0.00029 & 0.00020 & 0.00015 & 0.00003 & 0.00002 & 0.00000 & 0.00000 & 0.00634 \\
    32    & 0.00007 & 0.00005 & 0.00037 & 0.00028 & 0.00018 & 0.00013 & 0.00003 & 0.00002 & 0.00000 & 0.00000 & 0.00612 \\
    33    & 0.00007 & 0.00005 & 0.00035 & 0.00027 & 0.00016 & 0.00012 & 0.00002 & 0.00002 & 0.00000 & 0.00000 & 0.00592 \\
    151   & 0.00004 & 0.00003 & 0.00000 & 0.00000 & 0.00000 & 0.00000 & 0.00000 & 0.00000 & 0.00000 & 0.00000 & 0.00089 \\
    152   & 0.00004 & 0.00003 & 0.00000 & 0.00000 & 0.00000 & 0.00000 & 0.00000 & 0.00000 & 0.00000 & 0.00000 & 0.00088 \\
    153   & 0.00004 & 0.00003 & 0.00000 & 0.00000 & 0.00000 & 0.00000 & 0.00000 & 0.00000 & 0.00000 & 0.00000 & 0.00087 \\
    301   & 0.00001 & 0.00000 & 0.00000 & 0.00000 & 0.00000 & 0.00000 & 0.00000 & 0.00000 & 0.00000 & 0.00000 & 0.00011 \\
    302   & 0.00001 & 0.00000 & 0.00000 & 0.00000 & 0.00000 & 0.00000 & 0.00000 & 0.00000 & 0.00000 & 0.00000 & 0.00011 \\
    303   & 0.00001 & 0.00000 & 0.00000 & 0.00000 & 0.00000 & 0.00000 & 0.00000 & 0.00000 & 0.00000 & 0.00000 & 0.00011 \\
    501   & 0.00000 & 0.00000 & 0.00000 & 0.00000 & 0.00000 & 0.00000 & 0.00000 & 0.00000 & 0.00000 & 0.00000 & 0.00001 \\
    502   & 0.00000 & 0.00000 & 0.00000 & 0.00000 & 0.00000 & 0.00000 & 0.00000 & 0.00000 & 0.00000 & 0.00000 & 0.00001 \\
    503   & 0.00000 & 0.00000 & 0.00000 & 0.00000 & 0.00000 & 0.00000 & 0.00000 & 0.00000 & 0.00000 & 0.00000 & 0.00001 \\
  $\geq$504   & 0.00000 & 0.00000 & 0.00000 & 0.00000 & 0.00000 & 0.00000 & 0.00000 & 0.00000 & 0.00000 & 0.00000 & 0.00000 \\
    \hline
    Total & 0.0123 &	0.0094 &	0.02192 &	0.0164 &	0.0256 &	0.0192 &	0.0253  &	0.0190	& 0.02297 &	 0.0172 &	0.9999 \\
    \hline
   \end{tabular}%
  \label{maptab2}%
\end{table}%

\begin{table}[H]
\tiny
  \centering
  \caption{Joint probabilities (queue size and server content) at arbitrary epoch for SV}
  \smallskip
    \begin{tabular}{|p{0.8cm}|p{0.8cm}|p{0.8cm}|p{0.8cm}|p{0.8cm}|p{0.8cm}|p{0.8cm}|p{0.8cm}|p{0.8cm}|p{0.8cm}|p{0.8cm}|p{0.8cm}|p{0.8cm}|}
    \hline
    \smallskip
    $n$     & $R_{1}(n,0)$ & $R_{2}(n,0)$ & $\xi_{1}(n,5)$ & $\xi_{2}(n,5)$ & $\xi_{1}(n,6)$ & $\xi_{2}(n,6)$ & $\xi_{1}(n,7)$ & $\xi_{2}(n,7)$ & $\xi_{1}(n,8)$ & $\xi_{2}(n,8)$ & $\xi_{1}(n,9)$ & $\xi_{2}(n,9)$ \\
    \hline
    0     & 0.00000 & 0.00000 & 0.00421 & 0.00317 & 0.00268 & 0.00202 & 0.00220 & 0.00166 & 0.00179 & 0.00135 & 0.00324 & 0.00240 \\
    1     & 0.00002 & 0.00001 & 0.00346 & 0.00262 & 0.00209 & 0.00158 & 0.00162 & 0.00123 & 0.00125 & 0.00095 & 0.00383 & 0.00286 \\
    2     & 0.00011 & 0.00008 & 0.00273 & 0.00207 & 0.00154 & 0.00117 & 0.00113 & 0.00086 & 0.00082 & 0.00063 & 0.00400 & 0.00300 \\
    3     & 0.00046 & 0.00033 & 0.00208 & 0.00158 & 0.00110 & 0.00084 & 0.00075 & 0.00057 & 0.00051 & 0.00039 & 0.00392 & 0.00295 \\
    4     & 0.00126 & 0.00092 & 0.00154 & 0.00117 & 0.00076 & 0.00058 & 0.00049 & 0.00037 & 0.00031 & 0.00024 & 0.00373 & 0.00280 \\
    5     &       &       & 0.00112 & 0.00085 & 0.00051 & 0.00039 & 0.00031 & 0.00023 & 0.00018 & 0.00014 & 0.00348 & 0.00262 \\
    6     &       &       & 0.00079 & 0.00060 & 0.00034 & 0.00026 & 0.00019 & 0.00014 & 0.00011 & 0.00008 & 0.00323 & 0.00243 \\
    7     &       &       & 0.00056 & 0.00042 & 0.00022 & 0.00017 & 0.00011 & 0.00009 & 0.00006 & 0.00005 & 0.00299 & 0.00225 \\
    8     &       &       & 0.00038 & 0.00029 & 0.00014 & 0.00011 & 0.00007 & 0.00005 & 0.00003 & 0.00003 & 0.00277 & 0.00208 \\
    9     &       &       & 0.00026 & 0.00020 & 0.00009 & 0.00007 & 0.00004 & 0.00003 & 0.00002 & 0.00001 & 0.00257 & 0.00194 \\
    10    &       &       & 0.00018 & 0.00014 & 0.00006 & 0.00004 & 0.00002 & 0.00002 & 0.00001 & 0.00001 & 0.00240 & 0.00180 \\
    31    &       &       & 0.00000 & 0.00000 & 0.00000 & 0.00000 & 0.00000 & 0.00000 & 0.00000 & 0.00000 & 0.00099 & 0.00074 \\
    32    &       &       & 0.00000 & 0.00000 & 0.00000 & 0.00000 & 0.00000 & 0.00000 & 0.00000 & 0.00000 & 0.00096 & 0.00072 \\
    33    &       &       & 0.00000 & 0.00000 & 0.00000 & 0.00000 & 0.00000 & 0.00000 & 0.00000 & 0.00000 & 0.00093 & 0.00070 \\
    34    &       &       & 0.00000 & 0.00000 & 0.00000 & 0.00000 & 0.00000 & 0.00000 & 0.00000 & 0.00000 & 0.00091 & 0.00068 \\
    35    &       &       & 0.00000 & 0.00000 & 0.00000 & 0.00000 & 0.00000 & 0.00000 & 0.00000 & 0.00000 & 0.00089 & 0.00067 \\
    301   &       &       & 0.00000 & 0.00000 & 0.00000 & 0.00000 & 0.00000 & 0.00000 & 0.00000 & 0.00000 & 0.00002 & 0.00001 \\
    302   &       &       & 0.00000 & 0.00000 & 0.00000 & 0.00000 & 0.00000 & 0.00000 & 0.00000 & 0.00000 & 0.00002 & 0.00001 \\
    303   &       &       & 0.00000 & 0.00000 & 0.00000 & 0.00000 & 0.00000 & 0.00000 & 0.00000 & 0.00000 & 0.00002 & 0.00001 \\
    373   &       &       & 0.00000 & 0.00000 & 0.00000 & 0.00000 & 0.00000 & 0.00000 & 0.00000 & 0.00000 & 0.00001 & 0.00001 \\
    374   &       &       & 0.00000 & 0.00000 & 0.00000 & 0.00000 & 0.00000 & 0.00000 & 0.00000 & 0.00000 & 0.00001 & 0.00000 \\
    392   &       &       & 0.00000 & 0.00000 & 0.00000 & 0.00000 & 0.00000 & 0.00000 & 0.00000 & 0.00000 & 0.00001 & 0.00000 \\
    $\geq$ 393   &       &       & 0.00000 & 0.00000 & 0.00000 & 0.00000 & 0.00000 & 0.00000 & 0.00000 & 0.00000 & 0.00000 & 0.00000 \\
    \hline
    Total & 0.00185	& 0.00134	& 0.01767 &	0.01338 &	0.00962 &	0.00730 &	0.00696 &	0.00529 &	0.00511 &	 0.00389 & 0.12674 &	 0.09514\\
    \hline

    \end{tabular}%
  \label{maptab3}%
\end{table}%

\begin{table}[H]
\tiny
  \centering
  \caption{Joint probabilities (queue size and vacation type) at arbitrary epoch for SV}
  \smallskip
    \begin{tabular}{|p{0.8cm}|p{0.8cm}|p{0.8cm}|p{0.8cm}|p{0.8cm}|p{0.8cm}|p{0.8cm}|p{0.8cm}|p{0.8cm}|p{0.8cm}|p{0.8cm}|p{0.8cm}|}
    \hline
    \smallskip
    $n$     & $\gamma_{1}(n,0)$ & $\gamma_{2}(n,0)$ & $\gamma_{1}(n,1)$ & $\gamma_{2}(n,1)$ & $\gamma_{1}(n,2)$ & $\gamma_{2}(n,2)$ & $\gamma_{1}(n,3)$ & $\gamma_{2}(n,3)$ & $\gamma_{1}(n,4)$ & $\gamma_{2}(n,4)$ & $P_{n}^{queue}$ \\
    \hline
    0     & 0.00187 & 0.00132 & 0.00000 & 0.00000 & 0.00000 & 0.00000 & 0.00000 & 0.00000 & 0.00000 & 0.00000 & 0.02966 \\
    1     & 0.00185 & 0.00137 & 0.00323 & 0.00232 & 0.00000 & 0.00000 & 0.00000 & 0.00000 & 0.00000 & 0.00000 & 0.03354 \\
    2     & 0.00184 & 0.00138 & 0.00318 & 0.00237 & 0.00372 & 0.00269 & 0.00000 & 0.00000 & 0.00000 & 0.00000 & 0.03687 \\
    3     & 0.00184 & 0.00138 & 0.00314 & 0.00236 & 0.00357 & 0.00267 & 0.00356 & 0.00258 & 0.00000 & 0.00000 & 0.03983 \\
    4     & 0.00184 & 0.00138 & 0.00310 & 0.00232 & 0.00341 & 0.00256 & 0.00326 & 0.00244 & 0.00308 & 0.00224 & 0.04250 \\
    5     & 0.00183 & 0.00138 & 0.00304 & 0.00228 & 0.00323 & 0.00243 & 0.00292 & 0.00220 & 0.00262 & 0.00197 & 0.03589 \\
    6     & 0.00183 & 0.00137 & 0.00298 & 0.00224 & 0.00303 & 0.00228 & 0.00258 & 0.00195 & 0.00216 & 0.00163 & 0.03188 \\
    7     & 0.00183 & 0.00137 & 0.00292 & 0.00219 & 0.00283 & 0.00213 & 0.00225 & 0.00169 & 0.00173 & 0.00131 & 0.02841 \\
    8     & 0.00182 & 0.00137 & 0.00285 & 0.00214 & 0.00262 & 0.00197 & 0.00193 & 0.00146 & 0.00137 & 0.00103 & 0.02547 \\
    9     & 0.00182 & 0.00136 & 0.00277 & 0.00208 & 0.00242 & 0.00182 & 0.00165 & 0.00125 & 0.00106 & 0.00080 & 0.02299 \\
    10    & 0.00181 & 0.00136 & 0.00270 & 0.00203 & 0.00222 & 0.00167 & 0.00140 & 0.00105 & 0.00081 & 0.00062 & 0.02091 \\
    31    & 0.00164 & 0.00123 & 0.00118 & 0.00088 & 0.00023 & 0.00017 & 0.00002 & 0.00001 & 0.00000 & 0.00000 & 0.00714 \\
    32    & 0.00163 & 0.00122 & 0.00112 & 0.00084 & 0.00020 & 0.00015 & 0.00002 & 0.00001 & 0.00000 & 0.00000 & 0.00693 \\
    33    & 0.00162 & 0.00121 & 0.00107 & 0.00080 & 0.00018 & 0.00013 & 0.00001 & 0.00001 & 0.00000 & 0.00000 & 0.00672 \\
    34    & 0.00161 & 0.00121 & 0.00102 & 0.00077 & 0.00016 & 0.00012 & 0.00001 & 0.00001 & 0.00000 & 0.00000 & 0.00653 \\
    35    & 0.00160 & 0.00120 & 0.00097 & 0.00073 & 0.00014 & 0.00011 & 0.00001 & 0.00001 & 0.00000 & 0.00000 & 0.00635 \\
    301   & 0.00006 & 0.00004 & 0.00000 & 0.00000 & 0.00000 & 0.00000 & 0.00000 & 0.00000 & 0.00000 & 0.00000 & 0.00014 \\
    302   & 0.00006 & 0.00004 & 0.00000 & 0.00000 & 0.00000 & 0.00000 & 0.00000 & 0.00000 & 0.00000 & 0.00000 & 0.00013 \\
    303   & 0.00006 & 0.00004 & 0.00000 & 0.00000 & 0.00000 & 0.00000 & 0.00000 & 0.00000 & 0.00000 & 0.00000 & 0.00013 \\
    373   & 0.00002 & 0.00001 & 0.00000 & 0.00000 & 0.00000 & 0.00000 & 0.00000 & 0.00000 & 0.00000 & 0.00000 & 0.00005 \\
    374   & 0.00002 & 0.00001 & 0.00000 & 0.00000 & 0.00000 & 0.00000 & 0.00000 & 0.00000 & 0.00000 & 0.00000 & 0.00005 \\
    392   & 0.00001 & 0.00001 & 0.00000 & 0.00000 & 0.00000 & 0.00000 & 0.00000 & 0.00000 & 0.00000 & 0.00000 & 0.00003 \\
    393   & 0.00001 & 0.00001 & 0.00000 & 0.00000 & 0.00000 & 0.00000 & 0.00000 & 0.00000 & 0.00000 & 0.00000 & 0.00003 \\
    394   & 0.00001 & 0.00001 & 0.00000 & 0.00000 & 0.00000 & 0.00000 & 0.00000 & 0.00000 & 0.00000 & 0.00000 & 0.00003 \\
    446   & 0.00001 & 0.00000 & 0.00000 & 0.00000 & 0.00000 & 0.00000 & 0.00000 & 0.00000 & 0.00000 & 0.00000 & 0.00001 \\
    447   & 0.00001 & 0.00000 & 0.00000 & 0.00000 & 0.00000 & 0.00000 & 0.00000 & 0.00000 & 0.00000 & 0.00000 & 0.00001 \\
    514   & 0.00000 & 0.00000 & 0.00000 & 0.00000 & 0.00000 & 0.00000 & 0.00000 & 0.00000 & 0.00000 & 0.00000 & 0.00001 \\
    515   & 0.00000 & 0.00000 & 0.00000 & 0.00000 & 0.00000 & 0.00000 & 0.00000 & 0.00000 & 0.00000 & 0.00000 & 0.00001 \\
    $\geq 516$    & 0.00000 & 0.00000 & 0.00000 & 0.00000 & 0.00000 & 0.00000 & 0.00000 & 0.00000 & 0.00000 & 0.00000 & 0.00000 \\
    \hline
    Total & 0.20851	& 0.15635 &	0.09092 &	0.06817 &	0.04722 &	0.03541 &	0.02624 &	0.01969 &	0.01524 &	 0.01144 &	0.99994\\
    \hline

    \multicolumn{12}{c}{$L_q$=48.547,  $L_s$=51.133, $W_q$=0.859, $W_s$= 0.905, $L^{ser}$=8.141, $L^{vac}$=0.838
  }\\
  \multicolumn{12}{c}{ $P_{dor}$=0.0032, $P_{idle}$=0.682, $P_{busy}$=0.318}\\
\noalign{\smallskip}
    \end{tabular}%
  \label{maptab4}%
\end{table}%


\begin{table}[H]
  \tiny
  \centering
  \caption{Joint probabilities (queue size and server content) at service completion epoch for MV}
  \smallskip
    \begin{tabular}{|c|c|c|c|c|c|c|c|c|c|c|}
    \hline
    {$n$} & {$\xi_{1}^{+}(n,5)$} & {$\xi_{2}^{+}(n,5)$} & {$\xi_{1}^{+}(n,6)$} & {$\xi_{2}^{+}(n,6)$} & {$\xi_{1}^{+}(n,7)$} & {$\xi_{2}^{+}(n,7)$} & {$\xi_{1}^{+}(n,8)$} & {$\xi_{2}^{+}(n,8)$} & {$\xi_{1}^{+}(n,9)$} & {$\xi_{2}^{+}(n,9)$} \\
    \hline
    0     & 0.00188 & 0.00133 & 0.00233 & 0.00166 & 0.00264 & 0.00189 & 0.00283 & 0.00202 & 0.00293 & 0.00209 \\
    1     & 0.00332 & 0.00245 & 0.00381 & 0.00282 & 0.00402 & 0.00298 & 0.00402 & 0.00299 & 0.00647 & 0.00475 \\
    2     & 0.00393 & 0.00293 & 0.00416 & 0.00312 & 0.00408 & 0.00307 & 0.00381 & 0.00287 & 0.00916 & 0.00680 \\
    3     & 0.00387 & 0.00291 & 0.00379 & 0.00286 & 0.00345 & 0.00261 & 0.00301 & 0.00228 & 0.01066 & 0.00797 \\
    4     & 0.00343 & 0.00259 & 0.00310 & 0.00235 & 0.00262 & 0.00199 & 0.00214 & 0.00163 & 0.01118 & 0.00838 \\
    5     & 0.00283 & 0.00214 & 0.00237 & 0.00180 & 0.00186 & 0.00142 & 0.00142 & 0.00108 & 0.01103 & 0.00829 \\
    6     & 0.00223 & 0.00169 & 0.00172 & 0.00131 & 0.00126 & 0.00096 & 0.00089 & 0.00068 & 0.01052 & 0.00791 \\
    7     & 0.00169 & 0.00129 & 0.00121 & 0.00092 & 0.00082 & 0.00063 & 0.00054 & 0.00042 & 0.00985 & 0.00741 \\
    8     & 0.00125 & 0.00095 & 0.00083 & 0.00063 & 0.00052 & 0.00040 & 0.00032 & 0.00025 & 0.00914 & 0.00687 \\
    9     & 0.00090 & 0.00069 & 0.00055 & 0.00042 & 0.00032 & 0.00025 & 0.00019 & 0.00014 & 0.00845 & 0.00636 \\
    10    & 0.00064 & 0.00049 & 0.00036 & 0.00028 & 0.00020 & 0.00015 & 0.00011 & 0.00008 & 0.00782 & 0.00588 \\
    20    & 0.00001 & 0.00001 & 0.00000 & 0.00000 & 0.00000 & 0.00000 & 0.00000 & 0.00000 & 0.00431 & 0.00324 \\
    21    & 0.00001 & 0.00001 & 0.00000 & 0.00000 & 0.00000 & 0.00000 & 0.00000 & 0.00000 & 0.00413 & 0.00310 \\
    251   & 0.00000 & 0.00000 & 0.00000 & 0.00000 & 0.00000 & 0.00000 & 0.00000 & 0.00000 & 0.00012 & 0.00009 \\
    457   & 0.00000 & 0.00000 & 0.00000 & 0.00000 & 0.00000 & 0.00000 & 0.00000 & 0.00000 & 0.00001 & 0.00000 \\
    458   & 0.00000 & 0.00000 & 0.00000 & 0.00000 & 0.00000 & 0.00000 & 0.00000 & 0.00000 & 0.00001 & 0.00000 \\
    459   & 0.00000 & 0.00000 & 0.00000 & 0.00000 & 0.00000 & 0.00000 & 0.00000 & 0.00000 & 0.00001 & 0.00000 \\
    460   & 0.00000 & 0.00000 & 0.00000 & 0.00000 & 0.00000 & 0.00000 & 0.00000 & 0.00000 & 0.00000 & 0.00000 \\
    \hline
    Total & 0.02734 &	0.02051 &	0.02486 &	0.01865 &	0.02207 &	0.01656	& 0.01940 &	0.01455 &	0.36214 &	 0.27164\\
    \hline

    \end{tabular}%
  \label{maptab5}%
\end{table}%

\begin{table}[H]
\tiny
  \centering
  \caption{Joint probabilities (queue size and vacation type) at vacation termination epoch for MV}
  \smallskip
    \begin{tabular}{|p{0.8cm}|p{0.8cm}|p{0.8cm}|p{0.8cm}|p{0.8cm}|p{0.8cm}|p{0.8cm}|p{0.8cm}|p{0.8cm}|p{0.8cm}|p{0.8cm}|p{0.8cm}|}
    \hline
    {$n$} & {$\gamma_1^{+}(n,0)$} & {$\gamma_2^{+}(n,0)$} & {$\gamma_1^{+}(n,1)$} & {$\gamma_{2}^{+}(n,1)$} & {$\gamma_{1}^{+}(n,2)$} & {$\gamma_{2}^{+}(n,2)$} & {$\gamma_{1}^{+}(n,3)$} & {$\gamma_{2}^{+}(n,3)$} & {$\gamma_{1}^{+}(n,4)$} & {$\gamma_{2}^{+}(n,4)$} & {$P_{n}^{queue+}$} \\
    \hline
    0     & 0.00000 & 0.00000 & 0.00000 & 0.00000 & 0.00000 & 0.00000 & 0.00000 & 0.00000 & 0.00000 & 0.00000 & 0.02161 \\
    1     & 0.00001 & 0.00001 & 0.00009 & 0.00007 & 0.00000 & 0.00000 & 0.00000 & 0.00000 & 0.00000 & 0.00000 & 0.03780 \\
    2     & 0.00001 & 0.00001 & 0.00018 & 0.00013 & 0.00049 & 0.00035 & 0.00000 & 0.00000 & 0.00000 & 0.00000 & 0.04509 \\
    3     & 0.00001 & 0.00001 & 0.00025 & 0.00018 & 0.00084 & 0.00062 & 0.00133 & 0.00095 & 0.00000 & 0.00000 & 0.04760 \\
    4     & 0.00002 & 0.00001 & 0.00031 & 0.00023 & 0.00109 & 0.00081 & 0.00207 & 0.00153 & 0.00272 & 0.00194 & 0.05013 \\
    5     & 0.00002 & 0.00002 & 0.00036 & 0.00027 & 0.00125 & 0.00093 & 0.00242 & 0.00181 & 0.00375 & 0.00279 & 0.04785 \\
    6     & 0.00002 & 0.00002 & 0.00040 & 0.00030 & 0.00135 & 0.00101 & 0.00251 & 0.00188 & 0.00389 & 0.00292 & 0.04351 \\
    7     & 0.00003 & 0.00002 & 0.00044 & 0.00033 & 0.00140 & 0.00105 & 0.00245 & 0.00184 & 0.00360 & 0.00271 & 0.03862 \\
    8     & 0.00003 & 0.00002 & 0.00047 & 0.00035 & 0.00140 & 0.00105 & 0.00229 & 0.00172 & 0.00311 & 0.00235 & 0.03395 \\
    9     & 0.00003 & 0.00002 & 0.00049 & 0.00037 & 0.00138 & 0.00104 & 0.00208 & 0.00157 & 0.00259 & 0.00196 & 0.02979 \\
    10    & 0.00003 & 0.00003 & 0.00051 & 0.00038 & 0.00134 & 0.00101 & 0.00185 & 0.00139 & 0.00209 & 0.00158 & 0.02622 \\
    20    & 0.00006 & 0.00004 & 0.00052 & 0.00039 & 0.00065 & 0.00049 & 0.00034 & 0.00026 & 0.00013 & 0.00010 & 0.01054 \\
    21    & 0.00006 & 0.00004 & 0.00051 & 0.00038 & 0.00059 & 0.00044 & 0.00028 & 0.00021 & 0.00009 & 0.00007 & 0.00992 \\
    251   & 0.00001 & 0.00001 & 0.00000 & 0.00000 & 0.00000 & 0.00000 & 0.00000 & 0.00000 & 0.00000 & 0.00000 & 0.00023 \\
    457   & 0.00000 & 0.00000 & 0.00000 & 0.00000 & 0.00000 & 0.00000 & 0.00000 & 0.00000 & 0.00000 & 0.00000 & 0.00001 \\
    458   & 0.00000 & 0.00000 & 0.00000 & 0.00000 & 0.00000 & 0.00000 & 0.00000 & 0.00000 & 0.00000 & 0.00000 & 0.00001 \\
    459   & 0.00000 & 0.00000 & 0.00000 & 0.00000 & 0.00000 & 0.00000 & 0.00000 & 0.00000 & 0.00000 & 0.00000 & 0.00001 \\
    460   & 0.00000 & 0.00000 & 0.00000 & 0.00000 & 0.00000 & 0.00000 & 0.00000 & 0.00000 & 0.00000 & 0.00000 & 0.00001 \\
    500   & 0.00000 & 0.00000 & 0.00000 & 0.00000 & 0.00000 & 0.00000 & 0.00000 & 0.00000 & 0.00000 & 0.00000 & 0.00001 \\
    501   & 0.00000 & 0.00000 & 0.00000 & 0.00000 & 0.00000 & 0.00000 & 0.00000 & 0.00000 & 0.00000 & 0.00000 & 0.00001 \\
    502   & 0.00000 & 0.00000 & 0.00000 & 0.00000 & 0.00000 & 0.00000 & 0.00000 & 0.00000 & 0.00000 & 0.00000 & 0.00001 \\
    $\geq$ 503	 & 0.00000 &	0.00000	 & 0.00000 &	0.00000 &	0.00000 &	0.00000 &	0.00000 &	0.00000 &	 0.00000 &	0.00000 &	 0.00000\\
   \hline
   Total & 0.01235 &	0.00926 &	0.02160 &	0.01620 &	0.02577 &	0.01932 &	0.02720 &	0.02040 &	0.02864 &	 0.02148	& 0.99993 \\
   \hline
   \end{tabular}%
  \label{maptab6}%
\end{table}%

\begin{table}[H]
\tiny
  \centering
  \caption{Joint probabilities (queue size and server content) at arbitrary epoch for MV}
    \begin{tabular}{|c|c|c|c|c|c|c|c|c|c|c|}
   \hline
    $n$     &  $\xi_{1}(n,5)$ & $\xi_{2}(n,5)$ & $\xi_{1}(n,6)$ & $\xi_{2}(n,6)$ & $\xi_{1}(n,7)$ & $\xi_{2}(n,7)$ & $\xi_{1}(n,8)$ & $\xi_{2}(n,8)$ & $\xi_{1}(n,9)$ & $\xi_{2}(n,9)$ \\
    \hline
    0     & 0.00329 & 0.00247 & 0.00278 & 0.00209 & 0.00229 & 0.00173 & 0.00187 & 0.00141 & 0.00338 & 0.00250 \\
    1     & 0.00271 & 0.00205 & 0.00216 & 0.00164 & 0.00169 & 0.00128 & 0.00130 & 0.00099 & 0.00399 & 0.00298 \\
    2     & 0.00213 & 0.00162 & 0.00160 & 0.00122 & 0.00117 & 0.00089 & 0.00085 & 0.00065 & 0.00416 & 0.00312 \\
    3     & 0.00162 & 0.00123 & 0.00114 & 0.00087 & 0.00078 & 0.00060 & 0.00054 & 0.00041 & 0.00407 & 0.00306 \\
    4     & 0.00120 & 0.00091 & 0.00079 & 0.00060 & 0.00051 & 0.00039 & 0.00033 & 0.00025 & 0.00386 & 0.00290 \\
    5     & 0.00087 & 0.00066 & 0.00053 & 0.00041 & 0.00032 & 0.00024 & 0.00019 & 0.00015 & 0.00359 & 0.00270 \\
    6     & 0.00062 & 0.00047 & 0.00035 & 0.00027 & 0.00020 & 0.00015 & 0.00011 & 0.00009 & 0.00332 & 0.00250 \\
    7     & 0.00043 & 0.00033 & 0.00023 & 0.00017 & 0.00012 & 0.00009 & 0.00006 & 0.00005 & 0.00306 & 0.00230 \\
    8     & 0.00030 & 0.00023 & 0.00015 & 0.00011 & 0.00007 & 0.00005 & 0.00004 & 0.00003 & 0.00283 & 0.00213 \\
    9     & 0.00021 & 0.00016 & 0.00009 & 0.00007 & 0.00004 & 0.00003 & 0.00002 & 0.00001 & 0.00262 & 0.00197 \\
    10    & 0.00014 & 0.00011 & 0.00006 & 0.00004 & 0.00002 & 0.00002 & 0.00001 & 0.00001 & 0.00244 & 0.00183 \\
    101   & 0.00000 & 0.00000 & 0.00000 & 0.00000 & 0.00000 & 0.00000 & 0.00000 & 0.00000 & 0.00030 & 0.00023 \\
    102   & 0.00000 & 0.00000 & 0.00000 & 0.00000 & 0.00000 & 0.00000 & 0.00000 & 0.00000 & 0.00030 & 0.00022 \\
    389   & 0.00000 & 0.00000 & 0.00000 & 0.00000 & 0.00000 & 0.00000 & 0.00000 & 0.00000 & 0.00001 & 0.00000 \\
    390   & 0.00000 & 0.00000 & 0.00000 & 0.00000 & 0.00000 & 0.00000 & 0.00000 & 0.00000 & 0.00001 & 0.00000 \\
    391   & 0.00000 & 0.00000 & 0.00000 & 0.00000 & 0.00000 & 0.00000 & 0.00000 & 0.00000 & 0.00001 & 0.00000 \\
    $\geq$ 392   & 0.00000 & 0.00000 & 0.00000 & 0.00000 & 0.00000 & 0.00000 & 0.00000 & 0.00000 & 0.00000 & 0.00000 \\
    \hline
    Total & 0.01380	& 0.01045 &	0.00997 &	0.00756 &	0.00725 &	0.00551 &	0.00533 &	0.00406 &	0.12747 &	 0.09570\\
    \hline
  \end{tabular}%
  \label{maptab7}%
\end{table}%

\begin{table}[H]
\tiny
  \centering
  \caption{Joint probabilities (queue size and vacation type) at arbitrary epoch for MV}
    \begin{tabular}{|p{0.8cm}|p{0.8cm}|p{0.8cm}|p{0.8cm}|p{0.8cm}|p{0.8cm}|p{0.8cm}|p{0.8cm}|p{0.8cm}|p{0.8cm}|p{0.8cm}|p{0.8cm}|}
    \hline
    \smallskip
    $n$     & $\gamma_{1}(n,0)$ & $\gamma_{2}(n,0)$ & $\gamma_{1}(n,1)$ & $\gamma_{2}(n,1)$ & $\gamma_{1}(n,2)$ & $\gamma_{2}(n,2)$ & $\gamma_{1}(n,3)$ & $\gamma_{2}(n,3)$ & $\gamma_{1}(n,4)$ & $\gamma_{2}(n,4)$ & $P_{n}^{queue}$ \\
    \hline
    0     & 0.00185 & 0.00131 & 0.00000 & 0.00000 & 0.00000 & 0.00000 & 0.00000 & 0.00000 & 0.00000 & 0.00000 & 0.02851 \\
    1     & 0.00183 & 0.00136 & 0.00322 & 0.00231 & 0.00000 & 0.00000 & 0.00000 & 0.00000 & 0.00000 & 0.00000 & 0.03254 \\
    2     & 0.00183 & 0.00137 & 0.00317 & 0.00236 & 0.00378 & 0.00273 & 0.00000 & 0.00000 & 0.00000 & 0.00000 & 0.03601 \\
    3     & 0.00183 & 0.00137 & 0.00313 & 0.00234 & 0.00363 & 0.00271 & 0.00387 & 0.00280 & 0.00000 & 0.00000 & 0.03908 \\
    4     & 0.00182 & 0.00137 & 0.00308 & 0.00231 & 0.00346 & 0.00260 & 0.00353 & 0.00265 & 0.00388 & 0.00281 & 0.04183 \\
    5     & 0.00182 & 0.00136 & 0.00303 & 0.00227 & 0.00328 & 0.00246 & 0.00317 & 0.00239 & 0.00330 & 0.00248 & 0.03728 \\
    6     & 0.00182 & 0.00136 & 0.00297 & 0.00223 & 0.00308 & 0.00231 & 0.00280 & 0.00211 & 0.00272 & 0.00205 & 0.03311 \\
    7     & 0.00181 & 0.00136 & 0.00290 & 0.00218 & 0.00287 & 0.00216 & 0.00244 & 0.00184 & 0.00218 & 0.00165 & 0.02946 \\
    8     & 0.00181 & 0.00136 & 0.00283 & 0.00213 & 0.00266 & 0.00200 & 0.00210 & 0.00158 & 0.00172 & 0.00130 & 0.02635 \\
    9     & 0.00180 & 0.00135 & 0.00276 & 0.00207 & 0.00246 & 0.00185 & 0.00179 & 0.00135 & 0.00134 & 0.00101 & 0.02371 \\
    10    & 0.00180 & 0.00135 & 0.00268 & 0.00202 & 0.00226 & 0.00170 & 0.00151 & 0.00114 & 0.00102 & 0.00078 & 0.02149 \\
    101   & 0.00085 & 0.00064 & 0.00002 & 0.00002 & 0.00000 & 0.00000 & 0.00000 & 0.00000 & 0.00000 & 0.00000 & 0.00206 \\
    102   & 0.00084 & 0.00063 & 0.00002 & 0.00002 & 0.00000 & 0.00000 & 0.00000 & 0.00000 & 0.00000 & 0.00000 & 0.00204 \\
    389   & 0.00002 & 0.00001 & 0.00000 & 0.00000 & 0.00000 & 0.00000 & 0.00000 & 0.00000 & 0.00000 & 0.00000 & 0.00004 \\
    390   & 0.00001 & 0.00001 & 0.00000 & 0.00000 & 0.00000 & 0.00000 & 0.00000 & 0.00000 & 0.00000 & 0.00000 & 0.00004 \\
    391   & 0.00001 & 0.00001 & 0.00000 & 0.00000 & 0.00000 & 0.00000 & 0.00000 & 0.00000 & 0.00000 & 0.00000 & 0.00003 \\
    392   & 0.00001 & 0.00001 & 0.00000 & 0.00000 & 0.00000 & 0.00000 & 0.00000 & 0.00000 & 0.00000 & 0.00000 & 0.00003 \\
    458   & 0.00001 & 0.00000 & 0.00000 & 0.00000 & 0.00000 & 0.00000 & 0.00000 & 0.00000 & 0.00000 & 0.00000 & 0.00001 \\
    459   & 0.00001 & 0.00000 & 0.00000 & 0.00000 & 0.00000 & 0.00000 & 0.00000 & 0.00000 & 0.00000 & 0.00000 & 0.00001 \\
    460   & 0.00001 & 0.00000 & 0.00000 & 0.00000 & 0.00000 & 0.00000 & 0.00000 & 0.00000 & 0.00000 & 0.00000 & 0.00001 \\
    514   & 0.00000 & 0.00000 & 0.00000 & 0.00000 & 0.00000 & 0.00000 & 0.00000 & 0.00000 & 0.00000 & 0.00000 & 0.00001 \\
    515   & 0.00000 & 0.00000 & 0.00000 & 0.00000 & 0.00000 & 0.00000 & 0.00000 & 0.00000 & 0.00000 & 0.00000 & 0.00001 \\
    $\geq$ 516   & 0.00000 & 0.00000 & 0.00000 & 0.00000 & 0.00000 & 0.00000 & 0.00000 & 0.00000 & 0.00000 & 0.00000 & 0.00000 \\
    \hline
    Total & 0.20682	& 0.15507 &	0.09046 &	0.06782 &	0.04796 &	0.03596 &	0.02847 &	0.02135 &	0.01919 &	 0.01439 &	0.99994\\
    \hline
    \multicolumn{12}{c}{$L_q$=48.268,  $L_s$=50.837, $W_q$=0.854, $W_s$= 0.899, $L^{ser}$=8.224, $L^{vac}$=0.887
  }\\
  \multicolumn{12}{c}{ $P_{idle}$=0.687, $P_{busy}$=0.312}\\
\noalign{\smallskip}
\end{tabular}%
  \label{maptab8}%
\end{table}%

After tabular representation we present the behavior of the considered model in graphical form which represents the comparison of the queue size dependent vacation (QSDV) and queue size independent vacation (QSIV) policy. The following two cases are considered for the comparison scenario:\\
\textbf{Case 1.} The QSDV rates are considered as $\nu_k=(k+1)^{2}(1.1)$,~ $0\leq k\leq 2$.\\
\textbf{Case 2.} The QSIV rates are considered as $\nu_k=\nu_0$,~ $0\leq k\leq 2$. The vacation time decreases with increasing queue size at vacation initiation epoch in case 1. In case 2, however, no matter what the queue size is at vacation initiation, the vacation time remains constant. The following considerations apply to both cases (Case 1 and Case 2).
\begin{itemize}
	\item The vacation time of the server follows $E_2$ distribution, having PH-type representation ($\alpha_k, T_k$), $T_k$= $\begin{pmatrix}
	-\nu_k & \nu_k\\
	0.0 & -\nu_k
	\end{pmatrix}$,
	$\alpha_k$= $\begin{pmatrix}
	1.0 & 0.0
	\end{pmatrix}$, $0\leq k \leq 2$.

\item The service time of each batch under service follow the Erlang ($E_3$) distribution having PH-type representation ($\alpha_r, T_r$), where $T_r$= $\begin{pmatrix}
-\mu_r & \mu_r & 0.0\\
0.0 & -\mu_r & \mu_r\\
0.0 & 0.0 & -\mu_r
\end{pmatrix}$,  $\mu_r=0.3r$ , $\alpha_r$= $\begin{pmatrix}
1.0 & 0.0 & 0.0, \end{pmatrix}$, $3\leq r\leq 5$.
\item The $MAP$ is represented
by the matrices $C_l=$
$\begin{pmatrix}
-4.657l & 1.761l \\
1.128l & -3.941l
\end{pmatrix}$
and $D_l=$
$\begin{pmatrix}
1.657l & 1.239l\\
0.872l & 1.941l
\end{pmatrix}$. $\xi(C_l+D_l)$=$\textbf{0}$. $l$=1.0, 1.1,..., 2.0.  $\xi$=[0.4, 0.6].
\end{itemize}

\begin{figure}[H]
	\begin{minipage}[t]{5.0cm}
		\begin{center}
			\includegraphics[width=7cm]{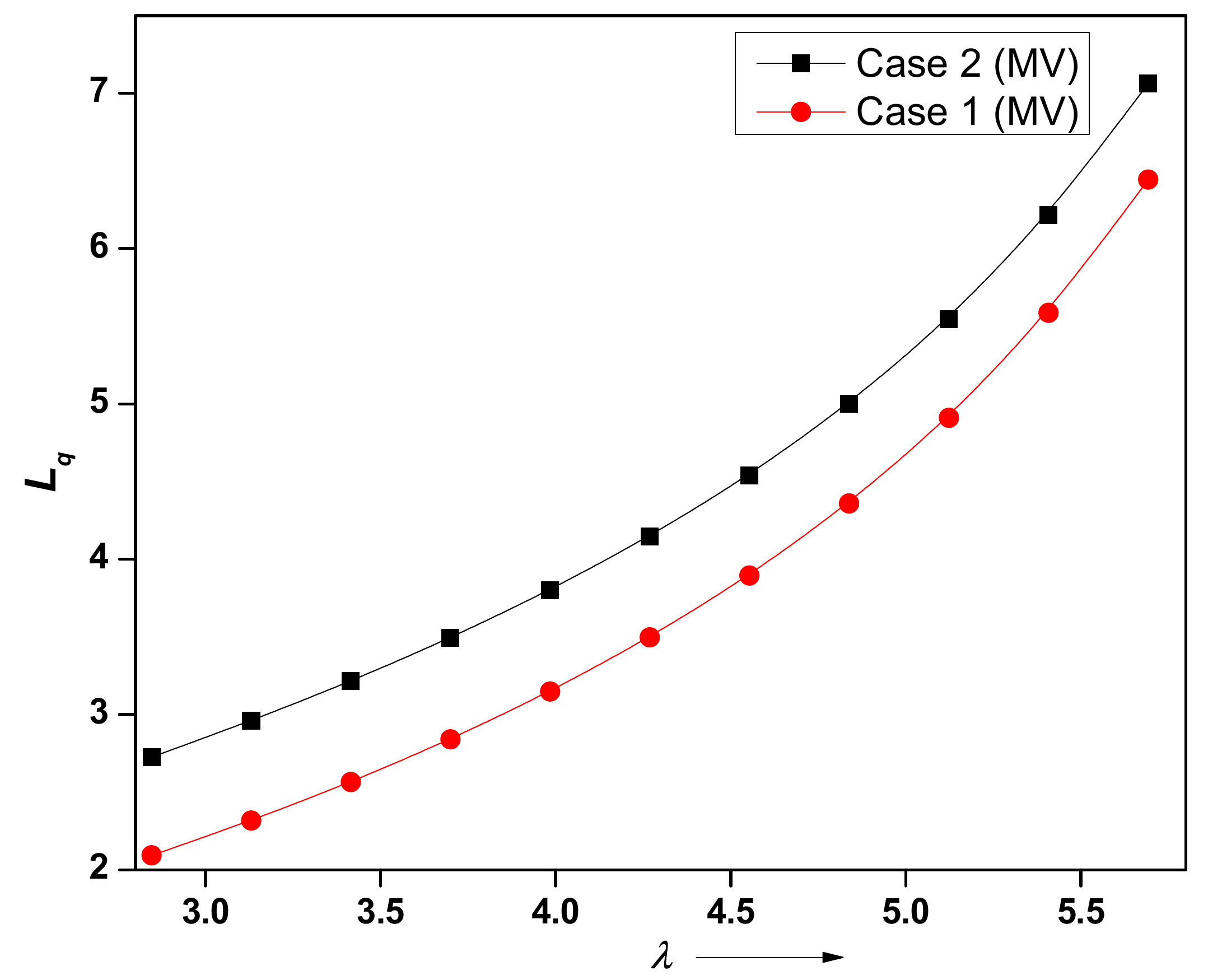}
			\vspace{-0.2 cm}{\caption[Short caption for figure 1]{\label{Fig1map} Effect of $\lambda$ on $L_q$}}
		\end{center}
	\end{minipage}
	\vspace{0.5 cm}
	\hfill
	\begin{minipage}[t]{5.0cm}
		\begin{center}
			\includegraphics[width=7cm]{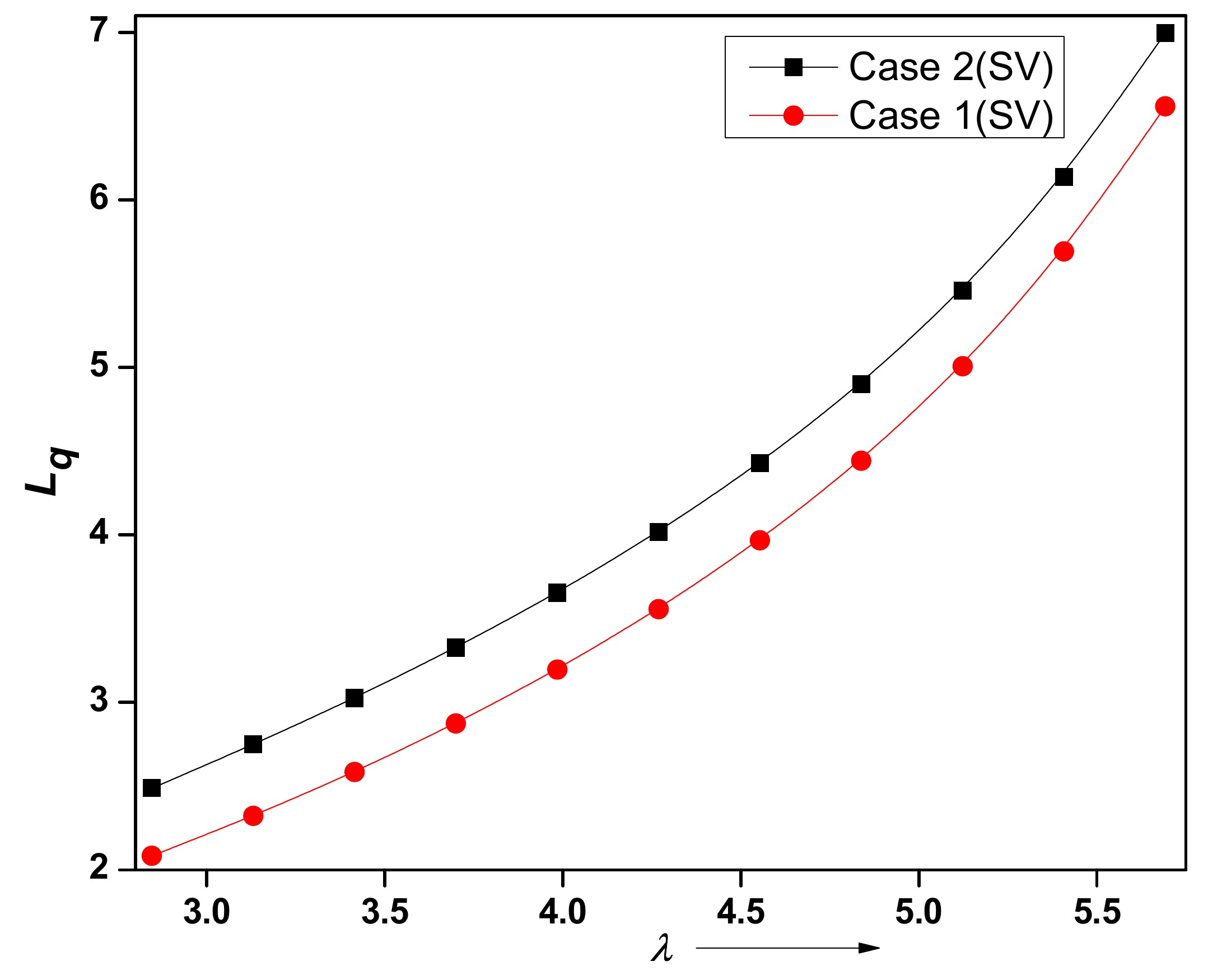}
			\vspace{-0.2 cm}{\caption[Short caption for figure 1]{\label{Fig2map} Effect of $\lambda$ on $L_q$}}
		\end{center}
	\end{minipage}
\hfill
\end{figure}

It is observed from Figure [\ref{Fig1map}-\ref{Fig2map}] that as the effective arrival rate $\lambda$ increases, the expected queue length $L_q$ increases in both the cases, this is because increasing the effective arrival rate increases the traffic intensity, and this behavior reflects the increase in $L_q$. Also, it can be marked here, $L_q$ is lower in Case 1 than Case 2 for a fixed $\lambda$.
Hence, the consideration of QSDV policy is more versed, because consideration of QSDV minimizes $L_q$ in comparison to the QSIV.

\section{Conclusion}\label{sec8} We discussed an infinite capacity $M/G_{r}^{(h,H)}/1$ queue with queue size dependent SV (MV). Bivariate vector generating function method and the supplementary variable approach have been used to extract steady state joint probabilities of queue content, server content (vacation type), and the phase of arrival process. The present model can be extended for analyzing different queueing systems with more different vacation policies (\textit{viz}., $BMAP/G_{r}^{(a,b)}/1$ queue with queue size dependent single and multiple working vacation), which is left for the future study.



\end{document}